\documentclass[10pt]{article}
\usepackage[margin=1in]{geometry} 
\usepackage{amsmath}        
\usepackage{graphicx}
\usepackage{amssymb}
\usepackage{amsthm}
\usepackage{epstopdf}
\usepackage{amsfonts}
\usepackage[scr]{rsfso}
\usepackage{hyperref}
\usepackage{booktabs}
\usepackage{placeins}
\makeatletter
\def\amsbb{\use@mathgroup \M@U \symAMSb}
\makeatother
\usepackage{bbold}
\usepackage{natbib}
\usepackage{epigraph}
\usepackage{authblk}
\usepackage{esint}
\usepackage{subfig}

\newtheorem{thm}{Theorem}[section]
\newtheorem{cor}{Corollary}[section]
\newtheorem{de}{Definition}[section]
\newtheorem{re}{Remark}[section]
\newtheorem{lem}{Lemma}[section]
\newtheorem{prop}{Proposition}[section]
\def\d{\delta}
\def\c{\circ}

\def\E{\amsbb{E}}

\def\C{\amsbb{C}ov}
\def\P{\amsbb{P}}
\def\R{\amsbb{R}}
\def\N{\amsbb{N}}

\def\D{\mathrm{d}}
\def\Var{\amsbb{V}ar}
\def\C{\amsbb{C}ov}

\newcommand{\norm}[1]{\left\lVert#1\right\rVert}
\providecommand{\norm}[1]{\lVert#1\rVert}
\DeclareMathOperator*{\argmax}{arg\,max}

\DeclareMathOperator*{\supp}{supp}

\DeclareMathOperator*{\for}{\quad\text{for}\quad}
\DeclareMathOperator*{\Leb}{Leb}

\DeclareMathOperator*{\diag}{diag}

\newcommand\ind[1]{\amsbb{I}_{#1}}

\usepackage{fancyhdr}
\usepackage{sectsty}
\usepackage{bold-extra}
\allsectionsfont{\bfseries\scshape}

\title{\scshape random graphs by product random measures}

\author[1]{Caleb Deen Bastian}
\author[1,2]{Herschel Rabitz}
\affil[1]{Program in Applied and Computational Mathematics, Princeton University, Princeton, NJ. 08544}
\affil[2]{Department of Chemistry, Princeton University, Princeton, NJ. 08544}

\date{\today}  

\begin{document}


\vspace{-6cm}                    
\maketitle

\begin{abstract}A natural representation of random graphs is the random measure. The collection of product random measures, their transformations, and non-negative test functions forms a general representation of the collection of non-negative weighted random graphs, directed or undirected, labeled or unlabeled, where (i) the composition of the test function and transformation is a non-negative edge weight function, (ii) the mean measures encode edge density/weight and vertex degree density/weight, and (iii) the mean edge weight, when square-integrable, encodes generalized spectral and Sobol representations. We develop a number of properties of these random graphs, and we give simple examples of some of their possible applications. 

\end{abstract}
\section{Introduction} 

Random graphs are ubiquitous. Following the pioneering work of \cite{er}, their forms have come to differentiate into a remarkable diversity of models, each characterized by its law or generative process. Well known classes of random graphs include namesake Erd\"{o}s-R\'{e}nyi, stochastic block \citep{block} and dot-product \citep{dot} models for social networks, the configuration model for fixed-degree graphs \citep{networks}, preferential attachment models such as Barb\'{a}si-Albert for generating scale-free graphs (those having power law degree distributions) \citep{pref}, Kallenberg exchangeable graphs invariant under arbitrary relabeling of the vertices, formed by exchangeable random measures on $\R_{\ge0}^2$ \citep{exchange3,exchange2}, exponential family models \citep{exp}, and so on. Random graphs also feature prominently in probabilistic graphical models, yielding the directed acyclic graphs of Bayesian networks and the undirected graphs of Gibbs (and Markov) networks \citep{pgm,gibbs}. In an attempt to make sense of the various models, taxonomies have been developed for organizing their concepts, such as their generators, metrics, and applications \citep{review}. 


To make a study of random graphs, we must have a definition for them. In the literature, there are multiple definitions. We represent the space of non-negative weighted $n$-graphs by the space of adjacency arrays indexed by $V^n$, where $V$ is an index set of vertices, and denote it as $H_n(V)\equiv\R_{\ge0}^{V^n}$. We use the convention that an edge possessing zero weight is regarded as `inactive' or absent. $V$ is either countable (finite or infinite) or uncountable (infinite). For each $n$ and $V$, the adjacency array $W\in H_n(V)$ encodes a non-negative weighted $n$-graph. Spaces can be combined through union, e.g. $H(V) \equiv \cup_{n\ge1} H_n(V)$ is the space of graphs containing edges of all orders, where $H_1(V)$ is the completely disconnected graph. Let $\pi$ be a probability measure on a space of adjacency arrays. Then a \emph{random graph} is identified to a random realization of an adjacency array having law $\pi$. One can either focus on defining $\pi$ explicitly or implicitly through the sampler, also known as a \emph{generator}. As a general strategy, we focus on generators. To comport with graph notation, we denote a random $n$-graph $G$ having adjacency array $W\in H_n(V)$ by the triple $G=(V,V^n,W)$.

There are a many ways of thinking about random graphs as objects. Recall a \emph{random field} is a random function or indexed collection of random variables. Thus a random adjacency array is a non-negative random field. Now consider integrals of the random adjacency array over the index dimensions. These integrals are \emph{random measures}. Such measures include edge count/weight and vertex degree count/weight. Thus a random graph is identified to a random field, and integrals over its adjacency array are random measures. Note that, conversely, random measures may be used to define random fields. These relations reflect reciprocities between random fields and random measures. 

Another perspective is to interpret a random field as a \emph{random transformation} of its domain or index set. This suggests a hierarchical model for random transformations, where their index sets may be random and labeled, with distributions specified for set cardinality and labels. Because random sets form random measures and vice versa, random graphs may be identified to random measures and their random transformations, where the top of the hierarchy is the random set, encoding (and encoded by) a random measure, with the next level of the hierarchy represented by the random transformation/field/function. This approach has the feature that the random measure is the base object, with the random transformation its secondary, so that integrals of the transformation are directly expressed in terms of the random measure, specifically its image, whose mean measure and Laplace functional respectively encode the integral means and laws.

Consider the non-negative weighted countable 2-graph $G=(V,V\times V,w)$ composed of a countable vertex set $V$, edge set $V\times V$, and non-negative weight function $w\ge0$, which forms the (adjacency array) $W=w(V\times V)=\{w(x,y): (x,y)\in V\times V\}$. Note that we have slightly generalized the graph notation, where the weight function is given instead of the adjacency array. As discussed, a general class of random graphs follows from considering a random product set $V\times V$ and a random function/field/transformation $w$. Special cases are deterministic $V\times V$ and/or $w$. The law of $w$ encodes the distribution of the edge weights conditioned on $V\times V$. The random set $V$ forms a \emph{random counting measure} $N$, also known as a \emph{point process}. Similarly, the random product set $V\times V$ forms a \emph{product random counting measure} $M=N\times N$. 

More generally, we denote the collection of non-negative weighted random 2-graphs as $G=G(w=g\circ f)=(V,V\times V,w=g\circ f)$ and compute its properties in terms of random measures $N$ on measurable space $(E,\mathscr{E})$ and $M$ on product space $(E\times E,\mathscr{E}\otimes\mathscr{E})$, whereby \begin{samepage}\begin{enumerate}\item[(i)] the graph edge set $V\times V$ is the Cartesian product of the vertex set $V$ (forming $N$) with itself, whose elements take values in $E\times E$ \item[(ii)]the edge weight function $w=g\circ f\ge0$ is the composition of a non-negative test function $g:F\mapsto\R_{\ge0}$ and (possibly random) transformation $f:E\times E\mapsto F$, where $(F,\mathscr{F})$ is a measurable space \item[(iii)] the mean measure of $M\circ f^{-1}$ encodes the mean edge count and weight of the graph model respectively as $\E(M\circ f^{-1})\ind{\supp(g)}$ and $\E(M\circ f^{-1})g$ \item[(iv)] putting $f_x(y)=f(x,y)$ and $f^x(y)=f(y,x)$, the mean measures of $N\circ f_x^{-1}$ and $N\circ f^{x-1}$ respectively encode the mean out- and in-degrees (counts) of vertex $x$ as $\E(N\circ f_x^{-1})\ind{\supp(g)}$ and $\E(N\circ f^{x-1})\ind{\supp(g)}$ and its mean weights as $\E(N\circ f_x^{-1})g$ and $\E(N\circ f^{x-1})g$ \item[(v)] the mean edge weight $\E w$, when square-integrable, encodes generalized spectral and Sobol representations of $G(w)$\end{enumerate}\end{samepage} 

These graphs may be directed or undirected, depending on the symmetry of $w$, and exist for general counting measures as well as for atomic and/or diffuse label distributions when labeled. The codomain of $w$ as subset of the non-negative reals determines graph type: $\{0,1\}$ corresponds to graphs, $\N_{\ge0}$ corresponds to multigraphs, and $\R_{\ge0}$ corresponds to weighted graphs proper. Extension to real weights forms signed measures. 

This representation is general: higher-order graphs and their superpositions may be provisioned using superpositions of higher-order product random measures. Uncountable graphs are retrieved by relaxing the countability of the vertex set $V$.  

In this context, of the aforementioned models, the Erd\"{o}s-R\'{e}nyi, stochastic block, dot-product, directed acyclic, and Kallenberg exchangeable random graphs can be formulated using random $V\times V$ and $w$, whereas fixed-degree, preferential attachment, exponential family, and Markov networks can be formulated using deterministic $V\times V$ and random $w$. 


The overall purpose of this article is to provide a perspective on thinking about random graphs as being formed by product random measures. Given the tremendous diversity of random graphs, we believe product random measures provide value as general representations equipped with calculus, serving as `scaffolding' for templating models and computing their properties. 




This article is organized as follows. Next, in a review section, we cover the mathematical methods and give key results on the random measures of interest (\ref{sec:methods}). Then we develop a general model of random graphs using product `stone-throwing' random measures and describe some graph properties (\ref{sec:rg}). We follow this with an exposition on random graphs formed by product random measures having fixed atoms and random non-negative (integer) weights (\ref{sec:farwg}). We discuss some simple applications to graphon identification (\ref{sec:graphons}), prime graphs (\ref{sec:primes}), spin networks (\ref{sec:spin}), Bayesian networks (\ref{sec:bayes}), and deep neural networks (\ref{sec:neural}). We end with discussion and conclusions (\ref{sec:discuss}).


\section{Review of methods / random measures} \label{sec:methods} Here we provide a review of the mathematical methods of random measures of this article. All the results are classical and well-known. For the reader who is knowledgeable on random measures, this section can be skipped (after suggested acquaintance with notation). For other readers, because these results will be our constant companions in later sections, we provide a reasonably thorough review of the calculus of random measures. First we lay out the backdrop of random measures (\ref{sec:background}). Then we give basic results on random measures based on the mixed binomial process called stone throwing construction (STC) random measures (\ref{sec:stc}), their traces (\ref{sec:trace}), and distributions (\ref{sec:dis}). We discuss random transformations (\ref{sec:rt}) and describe a product of a STC random measure with itself (\ref{sec:prod}).  Next we describe random measures having fixed atoms and non-negative random weights (FARW), their (counting) subclass having integer weights (FAIW), and their products (\ref{sec:farw}). We discuss a prototypical family of random counting measures called Poisson-type that are characterized by invariance under binomial thinning (\ref{sec:pt}). Lastly we briefly remark on random fields formed by random measures (\ref{sec:rf}). 

 
 \subsection{Backdrop}\label{sec:background}
 Let $(E,\mathscr{E})$ be a measurable space and let $(\Omega,\mathscr{H},\amsbb{P})$ be a probability space. A \emph{random measure} is a transition kernel from $(\Omega,\mathscr{H})$ into $(E,\mathscr{E})$. Specifically the mapping $N:\Omega\times E\mapsto\R_{\ge0}$ is a random measure if $\omega\mapsto N(\omega,A)$ is a random variable for each $A$ in $\mathscr{E}$ and if $A\mapsto N(\omega,A)$ is a measure on $(E,\mathscr{E})$ for each $\omega$ in $\Omega$. 
 
Let $\mathscr{E}_{\ge0}$ be the collection of non-negative $\mathscr{E}$-measurable functions. We give Fubini's foundational theorem of random measures. 
  
 \begin{samepage}
\begin{thm}[Fubini]\label{thm:fubini} Let $N$ be a random measure on $(E,\mathscr{E})$. Then \[Nf(\omega) = \int_E N(\omega,\D x)f(x)\for \omega\in\Omega\] defines a non-negative random variable $Nf$ for every function $f$ in $\mathscr{E}_{\ge0}$; \[\lambda(B)=\E N(B)=\int_\Omega\P(\D \omega)N(\omega,B)\for B\in\mathscr{E}\] defines a measure $\lambda=\E N$ on $(E,\mathscr{E})$ for each probability measure $\P$ on $(\Omega,\mathscr{H})$ called the mean or intensity measure;  and \[\E Nf = \lambda f=\int_\Omega\P(\D\omega)\int_E N(\omega,\D x)f(x)\] for every probability measure $\P$ on $(\Omega,\mathscr{H})$ and function $f$ in $\mathscr{E}_{\ge0}$.
\end{thm}
\end{samepage}\begin{proof} Theorem I.6.3 \citep{cinlar} \end{proof}
  
The law of $N$ is uniquely determined by the \emph{Laplace functional} $L:\mathscr{E}_{\ge0}\mapsto[0,1]$ \begin{equation}\label{eq:laplace} L(f)=\E e^{-Nf} = \E\exp_-\int_EN(\D x)f(x) \quad\text{for}\quad f\in\mathscr{E}_{\ge0}\end{equation} The Laplace functional encodes all the information of $N$: its distribution, moments, etc.  The distribution of $Nf$, i.e. $\P(Nf\in \D x)$, is encoded by the Laplace transform $\varphi$, which is expressed in terms of the Laplace functional \begin{equation}\varphi(\alpha) = \E e^{-\alpha Nf} =  \E e^{-N(\alpha f)}=L(\alpha f)\for \alpha\in\R_{\ge0}\end{equation} 


A random measure $N$ is said to be a \emph{random counting measure} if it is purely atomic and its every atom has weight one. Said another way, for each $A\subseteq E$, the random variable $N(A)$ is the number of points in $A$. A random counting measure $N$ is identified to a (countable) random point (multi)set $\mathbf{X}=\{X_i\}$, also known as the \emph{point process} associated with $N$. Some authors use random set notation, e.g., $\#(A)=N(A) = |\mathbf{X}\cap A|$. A variant of the random counting measure is allowing each atom to have a natural weight. We focus on two classes of purely atomic random counting measures: mixed binomial processes (MBP) having moving atoms conveyed through the stone throwing construction (STC) and random measures having fixed atoms and integer weights (FAIW). 


\subsection{Stone throwing}\label{sec:stc}
  
A broad class of purely atomic random counting measures are the \emph{mixed binomial processes} (MBP) \citep{kallenberg}. They are constructed as follows. Let $\nu$ be a probability measure on $(E,\mathscr{E})$, and let $\mathbf{X}=\{X_i\}$ be an independency (here a countable collection) of (iid) $E$-valued random variables having common law $\nu$, i.e. $X_i\sim\nu$ for the $i$. Let $\kappa$ be a $\N_{\ge 0}$-valued distribution with mean $c>0$ and variance $\d^2\ge0$, and let $K\sim\kappa$ be independent of $\mathbf{X}$. A MBP $N$ on $(E,\mathscr{E})$ is identified to a pair of deterministic probability measures $(\kappa,\nu)$ through the \emph{stone throwing construction}\citep{cinlar,Bastian:2020aa} (STC) \begin{equation}\label{eq:stc}N(A)=N\ind{A}=\int_{E}N(\D x)\ind{A}(x)\equiv\sum_{i}^{K}\ind{A}(X_i)\quad\text{for}\quad A\in\mathscr{E}\end{equation} where $\ind{A}$ is a set function. It is denoted the random measure $N=(\kappa,\nu)$ on $(E,\mathscr{E})$ and is said to be formed by $\mathbf{X}$. For functions $f$ in $\mathscr{E}_{\ge0}$, we have \[Nf = \int_EN(\D x)f(x)=\sum_i^Kf(X_i)\] Note that $\ind{A}(x)=\delta_x(A)$ so that $N$ be may concisely written as \[N = \sum_{i}^K\delta_{X_i}\] Below in Table~\ref{tab:dis} we give some frequently encountered counting measures. 


\begin{table}[h!]
\begin{center}

\begin{tabular}{lllc}
\toprule
Distribution & Parameters & Mass function (in $x$) & Support \\\midrule
Dirac & $n\in\N_{\ge1}$ & $\delta_n\{x\}$ & $\{n\}$ \\
Poisson & $c\in(0,\infty)$ & $c^xe^{-c}/x!$ & $\N_{\ge0}$\\
Negative binomial & $(r,p)\in\N_{\ge1}\times(0,1)$ & $\binom{r+x-1}{r-1}(1-p)^rp^x$  &$\N_{\ge0}$\\
Binomial & $(n,p)\in\N_{\ge1}\times(0,1)$ & $\binom{n}{x}(1-p)^{n-x}p^x$ & $\{0,1,\dotsb,n\}$\\
Uniform & $(m,n)\in\N_{\ge0}\times\N_{\ge1}$: $m\le n$& $1/(n-m+1)$ & $\{m,m+1,\dotsb,n-1,n\}$\\
Zeta & $s\in(1,\infty)$ & $x^{-s}/\zeta(s)$& $\N_{\ge1}$\\
Zipf & $(s,n)\in\R_{\ge0}\times\N_{\ge1}$ & $x^{-s}/H_x(s)$ & $\{1,\dotsb,n\}$\\
\bottomrule
\end{tabular}
\caption{Counting distributions $\kappa$}\label{tab:dis}
\end{center}
\end{table}

\FloatBarrier



The following theorem attains the Laplace functional of a STC random measure and the statistics of random variables formed by integrating test functions.

\begin{samepage}
\begin{thm}[STC]\label{thm:stc} Consider the counting random measure $N=(\kappa,\nu)$ on $(E,\mathscr{E})$ defined by STC \eqref{eq:stc}, where $K\sim\kappa$ has pgf $\psi$, mean $c>0$, and variance $\delta^2\ge0$. Then $N$ has Laplace functional \[L(f) = \E e^{- Nf} = \psi(\nu e^{-f})\for f\in\mathscr{E}_{\ge0}\] for function $f$ in $\mathscr{E}_{\ge0}$, $Nf$ has mean and variance \[\E Nf = c\nu f,\quad \Var Nf = c\nu f^2 + (\delta^2-c)(\nu f)^2\] and for functions $f, g$ in $\mathscr{E}_{\ge0}$, $Nf$ and $Ng$ have covariance  \[\C(Nf,Ng) = c\nu(fg) + (\delta^2-c)\nu f\,\nu g\]
\end{thm}\end{samepage}

\begin{proof} The Laplace functional follows from \[\E e^{-Nf} = \E(\E e^{-f(X)}) = \E(\nu e^{-f})^K = \psi(\nu e^{-f})\for f\in\mathscr{E}_{\ge0} \] The mean follows from Theorem~\ref{thm:fubini} with mean measure $c\nu$. The variance and covariance follow from the independence of $K$ and the $\{X_i\}$. 
\end{proof}

\subsection{Traces}\label{sec:trace}
Restrictions (traces) of random measures arise upon considering measurable subspaces. 

\begin{samepage}
\begin{thm}[Trace random measures]\label{thm:trace} Consider the random measure $N=(\kappa,\nu)$ on $(E,\mathscr{E})$ and the subspace $A\subseteq E$ with $\nu(A)=a>0$. Then \begin{enumerate}\item[(i)]$N_A(\cdot)=N(A\cap \cdot)=(N\ind{A},\nu_A)$ is the trace random measure of $N$ to $A$ on space $(E\cap A, \mathscr{E}_A)$, where $\nu_A(\cdot)=\nu(A\cap \cdot)/\nu(A)$ and $\mathscr{E}_A=\{A\cap B: B\in\mathscr{E}\}$\item[(ii)] $N_A$ has Laplace functional \[L_A(f) = \E e^{-N_A f} =\psi_A(\nu_A e^{-f})= \psi(a\nu_A e^{-f}+1-a)\quad\text{for}\quad f\in\mathscr{E}_{\ge0}\] where $\psi_A$ is pgf of $K_A=N\ind{A}$ and $\psi$ is pgf of $K=N\ind{E}$, and \item[(iii)] $N_Af$ has mean and variance \begin{align*}\E N_Af &=ac\nu_Af\\\Var N_A f &= ac\nu_A f^2 + a^2(\delta^2-c)(\nu_A f)^2\end{align*} and the covariance of $N_Af$ and $N_Ag$ is \[\C(N_Af,N_Ag) = ac\nu_A(fg) + a^2(\delta^2-c)\nu_A(f)\nu_A(g)\]
\end{enumerate} 

\end{thm}

\end{samepage}

\subsection{Distributions}\label{sec:dis} Next we give the distribution of the random variable $Nf$. The first relation is of course well known, following from the characteristic function, whereas the second on the Cauchy formula for extraction of coefficients from a whole series seems lesser so. 

\begin{samepage}
\begin{thm}[Distribution]\label{thm:dis} Let $N=(\kappa,\nu)$ be a random measure on $(E,\mathscr{E})$ and consider test function $f\in\mathscr{E}_{\ge0}$. Then the distribution of the random variable $Nf$ is given by \[\eta(x) = \frac{1}{2\pi}\int_\R \varphi(-i\alpha)e^{-i\alpha x}\D\alpha\for x\in\R_{\ge0}\] where $\varphi(\cdot)=L(\cdot f)=\psi(\nu e^{-\cdot f})$ is the Laplace transform of $Nf$. Now consider the indicator function $f=\ind{A}$, $A\subseteq E$ with $\nu(A)=a>0$. Then $K_A=Nf$ has pgf $\psi_A(t)=\psi(1-a+at)$ and distribution \[\eta\{x\}=\kappa_A\{x\}=\P(K_A=x) = \frac{1}{x!}\psi_A^{(x)}(0) = \frac{1}{2\pi i}\oint_{C}\frac{\psi_A(t)}{t^{x+1}}\D t\for x\in\N_{\ge0}\] where the contour integral is over the unit circle $C$.
\end{thm}
\end{samepage} \begin{proof} The law of $Nf$ is uniquely defined by the Laplace transform as \[\varphi(\alpha) = L(\alpha f)\for \alpha\in\R_{\ge0}\] where $L$ is the Laplace functional of $N$. It is also described through the characteristic function \[C(\alpha)=\varphi(-i\alpha)\for\alpha\in\R\] The density is attained through the inverse Fourier transform of the characteristic function. The pgf of $K_A$ is given by Theorem~\ref{thm:trace} as $\psi_A(t)=\psi(at+1-a)$. The contour integral formula follows from noting that \[\oint_C\frac{t^n}{t^{k+1}}\D t = \begin{cases}2\pi i &\text{if }n=k\\0 &\text{otherwise}\end{cases}\]
\end{proof}

We see that for indicator functions the probabilities are equivalently attained through derivatives and complex integration of the pgf, the latter being numerically far more stable. 

\subsection{Random transformations}\label{sec:rt} 

We define a random transformation. 

\begin{de}[Random transformation]\label{de:rt} Let $(E,\mathscr{E})$ and $(F,\mathscr{F})$ be measurable spaces. A \emph{random transformation} $\phi$ is a mapping from $E$ into $F$ \[\phi: (\omega,x)\mapsto\phi(\omega,x)\] that is measurable relative to $\mathscr{H}\otimes\mathscr{E}$ and $\mathscr{F}$. Writing $\phi x$ for the random variable $\omega\mapsto\phi(\omega,x)$, $\phi$ may be regarded as a collection of random variables $\phi=\{\phi x: x\in E\}$ with marginal distributions \[Q(x,B)=\P(\phi x\in B),\quad x\in E,\quad B\in\mathscr{F}\] that form, following from the joint measurability of $\phi$, a transition probability kernel $Q$ from $(E,\mathscr{E})$ into $(F,\mathscr{F})$. We refer to $Q$ as the marginal transition kernel of $\phi$. The law of $\phi$ is specified by the value $\P(\phi x\in A,\dotsb,\phi y\in B)$ for all finite collections of points $x,\dotsb,y$ in $E$ and sets $A,\dotsb,B$ in $\mathscr{F}$. If $\phi$ is an independency, then its law is uniquely determined by $Q$.
\end{de}

The following lemma is important: the mean measure of the image of a random measure under random transformation requires only the marginal transition kernel of the transformation, i.e. the transform marginal determines the mean measure, not the transform law. This of course is intuitive, as the actions of random measures are marginal integrals.

\begin{lem}[Mean]\label{lem:mean} Let $N$ be a random measure on $(E,\mathscr{E})$ with mean measure $\E N$. Let $\phi$ be a random transformation from $(E,\mathscr{E})$ into $(F,\mathscr{F})$ independent of $N$ with marginal transition kernel $Q$. Then the image of $N$ under $\phi$ has mean $\E(N\circ\phi^{-1})=(\E N)Q$.
\end{lem}
\begin{proof} The image random measure \[(N\circ\phi^{-1})f =\int_E N(\D x)f(\phi(x)) \for f\in\mathscr{F}_{\ge0}\] has mean \begin{align*}\E(N\circ\phi^{-1})f &=\E\int_E N(\D x)f(\phi(x))\\&=\int_E \E N(\D x)\E f(\phi(x))\\&=\int_{E}\E N(\D x)\int_FQ(x,\D z)f(z)\\&=(\E N)Qf\for f\in\mathscr{F}_{\ge0}\end{align*}
\end{proof}

On two dimensions the kernel's mean admits a generalized spectral representation. 

\begin{thm}[Generalized spectral representation]\label{re:spectral} Consider a random transformation $\phi:E\times E\mapsto\R_{\ge0}$ with mean function $f\in(\mathscr{E}\otimes\mathscr{E})_{\ge0}\cap  L^2(\nu\times\nu)$. Then \begin{enumerate}\item[(i)] $f$ defines a bounded, compact operator $T_f:L^2(\nu)\mapsto L^2(\nu)$ \[T_fg(x) = \int_Ef(x,y)g(y)\nu(\D y),\quad T_fg\in L^2(\nu)\for g\in L^2(\nu)\] \item[(ii)]$T_f$ admits a singular value decomposition \[T_f = \sum_n\sigma_n\langle\cdot,f_n\rangle g_n\] where $\langle\cdot,\cdot\rangle$ is the inner product of $L^2(\nu)$, $(\sigma_n)$ are non-negative (singular) values following $\lim_{n\rightarrow\infty}\sigma_n=0$, and $(f_n)$ and $(g_n)$ are orthonormal sets \item[(iii)] and \[f(x,y) = \sum_n\sigma_nf_n(x)g_n(y)\for (x,y)\in E\times E\]\end{enumerate} Moreover, if $f$ is symmetric, then $T_f$ is self-adjoint and $f_n=g_n$ for the $n$.
\end{thm} \begin{proof} By assumption, $f$ is a Hilbert-Schmidt kernel, so that $T_f$ is a Hilbert-Schmidt integral operator. Hence $T_f$ is continuous (and thus bounded) and compact, establishing (i). Because $T_f$ is a compact operator on the Hilbert space $L^2(\nu)$, it admits a singular value decomposition, establishing (ii) and (iii).\end{proof} 



\subsection{Product random measures of stone throwing}\label{sec:prod}
Often it is incumbent to think of the random measure $M=N\times N$ formed by the product of a random measure $N$ with itself, which we refer to as the product random measure $M$ of $N$. The following theorem calculates the mean (intensity) measure of the product random measure $M$ of a general STC $N$, as well as its image under transformation. The mean measure of $M$ is the non-factorial second moment measure of $N$. 

\begin{samepage}
\begin{thm}[Product STC]\label{thm:product} Consider the random measure STC $N=(\kappa,\nu)$ on $(E,\mathscr{E})$ formed from the independency $\{X_i:i=1,\dotsb,K\}$ where $K\sim\kappa$ with mean $c$ and variance $\delta^2$. Form the product random measure $M=N\times N$ on $(E\times E,\mathscr{E}\otimes\mathscr{E})$ as \begin{equation}\label{eq:product}Mf = \int_{E\times E}N(\D x)N(\D y)f(x,y) = \sum_i^K\sum_j^Kf\circ(X_i,X_j)\for f\in(\mathscr{E}\otimes\mathscr{E})_{\ge0}\end{equation}Then $M$ has mean \[\E M f = c(\nu\times I)f+(c^2+\delta^2-c)(\nu\times\nu)f\for f\in(\mathscr{E}\otimes\mathscr{E})_{\ge0}\] where $I$ is the identity kernel $I(x,\cdot)=\delta_x(\cdot)$. Moreover, for measurable space $(F,\mathscr{F})$ and random transformation $\phi:E\times E\mapsto F$ independent of $M$ with marginal transition kernel $Q$, the image random measure $M\circ\phi^{-1}$ on $(F,\mathscr{F})$ \begin{equation}\label{eq:image}(M\circ\phi^{-1})f = \int_{E\times E}N(\D x)N(\D y)f(\phi(x,y)) = \sum_i^K\sum_j^Kf\circ\phi(X_i,X_j)\for f\in\mathscr{F}_{\ge0}\end{equation} has mean \[\E(M\circ\phi^{-1})f = c(\nu\times I)Qf + (c^2+\delta^2-c)(\nu\times\nu)Qf\for f\in\mathscr{F}_{\ge0}\]
\end{thm}
\end{samepage} 

 \begin{proof}

Note that \[\E NfNg - \E N f \E N g=\C(Nf,Ng)=c\nu(fg)+(\delta^2-c)\nu(f)\nu(g)\] so that \[\E NfNg = c\nu(fg)+(c^2+\delta^2-c)\nu(f)\nu(g)\] Now take $f=\ind{A}$ and $g=\ind{B}$, giving \[\E M(A\times B) = c\nu(A\cap B) +(c^2+\delta^2-c)(\nu\times\nu)(A\times B)\] Taking $f=\ind{A\times B}\in(\mathscr{E}\otimes\mathscr{E})_{\ge0}$ and applying a monotone class argument gives \[\E M f = c(\nu\times I)f+(c^2+\delta^2-c)(\nu\times\nu)f\for f\in(\mathscr{E}\otimes\mathscr{E})_{\ge0}\] The mean of the image measure follows from Lemma~\ref{lem:mean} as $\E(M\circ\phi^{-1})=(\E M)Q$.

\end{proof}

We give a series of remarks on the theorem, the first on notation.

\begin{re}[Notation]The expressions of Theorem~\ref{thm:product} are given by \[\E Mf = c\int_E\nu(\D x)f(x,x)+(c^2+\delta^2-c)\int_{E\times E}\nu(\D x)\nu(\D y)f(x,y)\for f\in(\mathscr{E}\otimes\mathscr{E})_{\ge0}\] and  \[\E (M\circ\phi^{-1})f = c\int_{E\times F}\nu(\D x)Q((x,x),\D z)f(z) + (c^2+\delta^2-c)\int_{E\times E\times F}\nu(\D x)\nu(\D y)Q((x,y),\D z)f(z)\for f\in\mathscr{F}_{\ge0}\]
\end{re}

Next we give the mean measures and discuss deterministic transformations.

\begin{re}[Mean measures] The mean measure of $M$ is \[\E M = c(\nu\times I)+(c^2+\delta^2-c)(\nu\times\nu)\] The mean measure of $M\circ\phi^{-1}$ is \[\E(M\circ\phi^{-1}) = c(\nu\times I)Q + (c^2+\delta^2-c)(\nu\times\nu)Q\]
\end{re} 

\begin{re}[Deterministic $\phi$] The deterministic transformation $\phi:E\times E\mapsto F$ has marginal transition kernel $Q((x,y),\cdot)=\delta_{\phi(x,y)}(\cdot)$, and the mean measure is given by \[\E(M\circ\phi^{-1}) = (\E M)\circ\phi^{-1} = c(\nu\times I)\circ\phi^{-1} + (c^2+\delta^2-c)(\nu\times\nu)\circ\phi^{-1}\] 
\end{re}

We remark that the mean measure $\E M$ defines the second moment measures of $N$. 

\begin{re}[Second moment measures]\label{re:second} $\E M$ is the non-factorial  second moment measure of $N$. $\E\widetilde{M}=\E M - c(\nu\times I)=(c^2+\delta^2-c)(\nu\times\nu)$ is the factorial second moment measure of $N$, where $\E K(K-1)=c^2+\delta^2-c$ is the factorial second moment of $K$. 
\end{re}



A prototypical random measure is Poisson. 

\begin{re}[Poisson] For Poisson, put $\lambda = c\,\nu$. Then $\E M = \lambda\times I + \lambda\times\lambda$. This result holds more generally for infinite measures $\lambda$ on $(E,\mathscr{E})$ (in the infinite setting such measures on $(\R^n,\mathscr{B}_{\R^N})$ are referred to as \emph{L\'{e}vy measures}, where $\lambda f<\infty$ holds for $f(x)=|x|\wedge 1$ with $a\wedge b \equiv \min\{a,b\}$). For example, consider the Poisson random measure $N$ on $(E,\mathscr{E})=(\R_{\ge0},\mathscr{B}_{\R_{\ge0}})$ with infinite mean measure $\lambda=c\Leb$, $c>0$. Consider arbitrary $A,B\subseteq E$. Then $\E M(A\times B) = c|A\cap B| +c^2|A||B|$ and $\E \widetilde{M}(A\times B)=c^2|A||B|$ are respectively the non-factorial and factorial second moment measures of $N$. 
\end{re}

An important (sub)collection of functions are separable. For these, we attain the random measure variance.
\begin{re}[Separable functions] The subcollection of separable functions is denoted $\mathscr{E}_{\ge0}^2=\mathscr{E}_{\ge0}\times\mathscr{E}_{\ge0}\subset(\mathscr{E}\otimes\mathscr{E})_{\ge0}$ and defined such that every separable $h\in\mathscr{E}_{\ge0}^2$ is represented as $h=f\times g$ for $f,g\in\mathscr{E}_{\ge0}$, i.e. $h(x,y)=f(x)g(y)$. 
\end{re}

\begin{re}[Variance of separable functions] Consider a separable function $h=f\times g$. Then the moments of $Mh=(N\times N)(f\times g)=NfNg$ may be obtained through the Laplace functional of $N$, i.e. the first moment is \begin{align*}\E Mh &= \lim_{q,r\rightarrow0}\frac{\partial}{\partial q}\frac{\partial}{\partial r}\E e^{-N(qf+rg)} \end{align*} which by Theorem~\ref{thm:product} is explicitly given as \[\E Mh = c\nu(fg)+(c^2+\delta^2-c)\nu(f)\nu(g)\] and the second moment is \[\E (Mh)^2 = \lim_{q,r\rightarrow0}\frac{\partial^2}{\partial q^2}\frac{\partial^2}{\partial r^2}\E e^{-N(qf+rg)}\] which first and second moments define the variance $\Var(Mh) = \E (Mh)^2 - (\E Mh)^2$.
\end{re}

We give a remark on the normalized product random measure.

\begin{re}[Normalized product STC]\label{re:norm} The normalized product random measure $M^*$ of $M$ is defined as \[M^*f = \sum_i^Kf(X_i,X_i) + \frac{1}{2}\sum_{i\ne j}^Kf(X_i,X_j)\for f\in(\mathscr{E}\otimes\mathscr{E})_{\ge0}\] For random transformation $\phi$ with marginal transition kernel $Q$, the mean of the image $M^*\circ\phi^{-1}$ is  \[\E(M^*\circ\phi^{-1})f = c(\nu\times I)Qf + \frac{1}{2}(c^2+\delta^2-c)(\nu\times\nu)Qf\for f\in\mathscr{F}_{\ge0}\]
\end{re}

%

\subsection{Product random measures having fixed atoms and random weights}\label{sec:farw} 


Here we describe the class of random measures having fixed atoms and random weights (FARW). First we give its definition and basic results in the following theorem and give a few immediate remarks.

\begin{thm}[FARW]Let $D\subseteq E$ be a countable subset of $E$, and let $\{W_x: x\in D\}$ be an independency of non-negative random variables distributed $W_x\sim\nu_x$ with mean $c_x$ and variance $\delta_x^2$. Consider the random measure $N$ on $(E,\mathscr{E})$ formed as \begin{equation}\label{eq:fixed}N(A) = \int_{E}N(\D x)\ind{A}(x) = \sum_{x\in D}W_x\ind{A}(x)\for A\in\mathscr{E}\end{equation} with \emph{fixed atoms} of $D$ and \emph{random weights} of $\{W_x\}$. Then the Laplace functional of $N$ is given by \[L(f)=\E e^{-Nf} = \prod_{x\in D}\varphi_x(f(x))\for f\in\mathscr{E}_{\ge0}\] where $\varphi_x$ is the Laplace transform of $\nu_x$; for the function $f$ in $\mathscr{E}_{\ge0}$, $Nf$ has mean and variance \[\E Nf = \sum_{x\in D}c_xf(x),\quad \Var Nf = \sum_{x\in D}\delta_x^2f^2(x)\] and for the functions $f,g$ in $\mathscr{E}_{\ge0}$, $Nf$ and $Ng$ have covariance \[\C(Nf,Ng) = \sum_{x\in D}\delta_x^2f(x)g(x)\] 
\end{thm}

\begin{re}[Finiteness]If $\E Nf <\infty$, then $Nf<\infty$ almost surely.
\end{re}

\begin{re}[Poisson] Consider $\nu_x=\text{Poisson}(c_x)$. Then the Laplace functional \[L(f) =\E e^{-Nf} = \exp_-\sum_{x\in D}c_x(1-e^{-f(x)})=\exp_-\lambda(1-e^{-f})\for f\in\mathscr{E}_{\ge0}\] shows $N$ is a Poisson random measure on $(E,\mathscr{E})$ with mean measure $\lambda = \sum_{x\in D}c_x\delta_{x}$. When $\lambda$ is finite, $N$ is finite.
\end{re}

\begin{re}[Restriction]Consider subset $A\subseteq E$. Then the restricted random measure $N_A(\cdot)=N(A\cap\cdot)$ has Laplace functional \[L_A(f) = \E e^{-N_Af} = \prod_{x\in D\cap A}\varphi_x(f(x))\for f\in\mathscr{E}_{\ge0}\]
\end{re}

An important subclass of FARW random measures are those with integer weights. 
\begin{re}[Counting/FAIW]\label{re:count} The FARW random measure $M$ on $(E,\mathscr{E})$ is a counting measure when the weights are non-negative integers (FAIW): $W_x\sim\kappa_x$ for counting measure $\kappa_x$.
\end{re}

 We have the following theorem and remark on the mean measure of the product of a FARW random measure with itself. 
 
\begin{samepage}
\begin{thm}[Product FARW]\label{thm:product2} Consider the random measure $N$ on $(E,\mathscr{E})$ with countable set of fixed atoms $D\subseteq E$ and independency of non-negative random weights $\{W_x: x\in D\}$ distributed $W_x\sim\nu_x$ with mean $c_x$ and variance $\delta_x^2$. Form the product random measure $M=N\times N$ on $(E\times E,\mathscr{E}\otimes\mathscr{E})$ as  \begin{equation}\label{eq:product2}Mf =\int_{E\times E}N(\D x)N(\D y)f(x,y)= \sum_{(x,y)\in D\times D} W_xW_yf(x,y)\for f\in(\mathscr{E}\otimes\mathscr{E})_{\ge0}\end{equation}Then $M$ has mean \[\E Mf =\sum_{(x,y)\in D\times D}Z_{xy}f(x,y)\for f\in(\mathscr{E}\otimes\mathscr{E})_{\ge0}\] where \[Z_{xy}=\begin{cases}c_x^2+\delta_x^2&\text{if }x=y\\c_xc_y &\text{otherwise}\end{cases}\] Moreover, for measurable space $(F,\mathscr{F})$ and random transformation $\phi:E\times E\mapsto F$ independent of $M$ with marginal transition kernel $Q$, putting $Q_{xy}(\cdot)=Q((x,y),\cdot)$ the image random measure $M\circ\phi^{-1}$ on $(F,\mathscr{F})$ \begin{equation}\label{eq:image}(M\circ\phi^{-1})f = \int_{E\times E}N(\D x)N(\D y)f(\phi(x,y)) = \sum_{(x,y)\in D\times D} W_xW_yf\circ\phi(x,y)\for f\in\mathscr{F}_{\ge0}\end{equation} has mean \[\E(M\circ\phi^{-1})f = \sum_{(x,y)\in D\times D}Z_{xy}Q_{xy}(f)\for f\in\mathscr{F}_{\ge0}\]
\end{thm}
\end{samepage} 

\begin{re}[Mean measures] The mean measures of $M$ and $M\circ\phi^{-1}$ are respectively \[\E M = \sum_{(x,y)\in D\times D}Z_{xy}\] and \[\E(M\circ\phi^{-1}) = \sum_{(x,y)\in D\times D}Z_{xy}Q_{xy}\]
\end{re}

The FARW random measure $M$ encodes a random transform of a non-negative double array $B=f\circ\phi(D\times D)$. 

\begin{de}[$W$-transform]\label{re:wtransform} Consider the set-up and conclusions of Theorem~\ref{thm:product2}. Form the non-negative double array $B\equiv f\circ \phi(D\times D)=\{B_{xy}=f\circ\phi(x,y): (x,y)\in D\times D\}$ and the non-negative row array $W=\{W_x: x\in D\}$ so that (in matrix notation) $(M\circ\phi^{-1})f = W B W^\intercal=\norm{\diag(W)B\diag(W)}_1$ where $\norm{\cdot}_1$ is the entrywise 1-norm and $\diag(W)\equiv\{W_x\ind{}(x=y): (x,y)\in D\times D\}$. We refer to \[B_W=T_W(B)\equiv\diag(W)B\diag(W)\equiv\{W_xW_yB_{xy}: (x,y)\in D\times D\}\] as the $W$-transform $T_W$ of $B$ with mean $\E B_W =\E T_W(B)\equiv \{Z_{xy}Q_{xy}(f):(x,y)\in D\times D\}$. $T_W$ rescales the rows and columns of its double array argument relative to $W$. $T_W$ expresses the identity and zero transformations: for $W^1\equiv\{1: x\in D\}$, then $B_{W^1}=T_{W^1}(B)=B$; for $W^0\equiv\{0:x\in D\}$, then $B_{W^0}=T_{W^0}(B)=\{0:(x,y)\in D\times D\}$. When $D$ has finite size $n$, then $B$ and $B_W$ are $n\times n$ square matrices.  
\end{de}

The Laplace functional of $M$ may be attained. 

\begin{thm}[Product FARW Laplace functional]\label{prop:lf} Let $N$ be a FARW random measure on $(E,\mathscr{E})$ with countable fixed atoms $D\subseteq E$ and independent non-negative random weights $W=\{W_x:x\in D\}$ in $H=\R_{\ge0}^{|D|}$ distributed $W_x\sim\nu_x$. Then the Laplace functional of $M=N\times N$ is given by \[L(f) = \E e^{-Mf} = \int_H\left(\prod_{x\in D}\nu_x(\D w_x)\right)e^{-\sum_{(x,y)\in D\times D}w_xw_yf(x,y)}\for f\in(\mathscr{E}\otimes\mathscr{E})_{\ge0}\] 
\end{thm}
\begin{proof} By definition, \begin{align*}\E e^{-Mf} &= \E e^{-\sum_{(x,y)\in D\times D}W_xW_yf(x,y)}\\&=\int_H\nu(\D w) e^{-\sum_{(x,y)\in D\times D}w_xw_yf(x,y)}\\&=\int_H\left(\prod_{x\in D}\nu_x(\D w_x)\right) e^{-\sum_{(x,y)\in D\times D}w_xw_yf(x,y)}\for f\in(\mathscr{E}\otimes\mathscr{E})_{\ge0}\end{align*}
\end{proof}

This establishes an equivalence between product FARW random measures and Gibbs random fields.

\begin{re}[Partition function/Gibbs measure]\label{re:part} The Laplace transform $\varphi$ of the product FARW random measure $M$, i.e. $\varphi(\beta)=L(\beta f)$, $\beta\ge0$, is identical to the partition function $\mathcal{Z}$ of a Gibbs family of (probability) measures $(\P_\beta)_{\beta\ge0}$ with energy function $E_W\equiv Mf$ i.e. \[\P_\beta(W=w) = \frac{1}{\mathcal{Z}(\beta)}e^{-\beta E_w}=\frac{1}{\varphi(\beta)}e^{-\beta E_w}\for\beta\in\R_{\ge0}\] Hence the models are stochastically equivalent.
\end{re}

We give an example of a finite product FAIW random measure with countable atoms. 
\begin{re}[Product FAIW example]\label{re:finite} Consider product FAIW random measure $M=N\times N$ on $(E\times E,\mathscr{E}\otimes\mathscr{E})$ with Bernoulli weights (of $N$) distributed $W_x\sim\kappa_x=\text{Bernoulli}(p_x)$ such that \[Z_{xy}=\E W_xW_y = \begin{cases}p_x&\text{if }x=y\\p_xp_y&\text{otherwise}\end{cases}\]  Take $E=D=\N_{\ge1}$. $M$ is a counting measure. For each $x\in D$, let $p_x=1/x^s$, $s>1$. Then \[\E M(E\times E)= \zeta(s)-\zeta(2s)+\zeta^2(s) <\infty\] where $\zeta$ is the zeta function. Hence for this parameterization of $p_x$, $M$ is finite almost surely.
\end{re}

We have a slight generalization of the product random measure. 
\begin{prop}[Product FARW generalization]\label{prop:gen}Let $N$ be a FARW random measure. Then the product random measure $M$ defined as \[M(f,g) = (N\times N)f + Ng \for f\in(\mathscr{E}\otimes\mathscr{E})_{\ge0},\quad g\in\mathscr{E}_{\ge0}\] has Laplace functional \[L(f,g) = \E e^{-M(f,g)} =\int_H\left(\prod_{x\in D}\nu_x(\D w_x)\right) e^{-\left(\sum_{(x,y)\in D\times D}w_xw_yf(x,y) + \sum_{x\in D}w_xg(x)\right)}\for f\in(\mathscr{E}\otimes\mathscr{E})_{\ge0},\quad g\in\mathscr{E}_{\ge0}\]
\end{prop}

As with product STC random measures, we can define a normalized product FARW random measure.

\begin{re}[Normalized product FARW]\label{re:norm2} Let $M=N\times N$ be a product FARW random measure. The normalized product random measure $M^*$ of $M$ is defined as \[M^*f = \sum_{x}W^2_xf(x,x) + \frac{1}{2}\sum_{x\ne y}W_xW_yf(x,y)\for f\in(\mathscr{E}\otimes\mathscr{E})_{\ge0}\] For random transformation $\phi$ with marginal transition kernel $Q$, the image of $M^*$ under $\phi$ has mean \[\E (M^*\circ\phi^{-1})f = \sum_{x}Z_{xx}Q_{xx}(f) + \frac{1}{2}\sum_{x\ne y}Z_{xy}Q_{xy}(f)\for f\in\mathscr{F}_{\ge0}\]
\end{re}



\subsection{Poisson-type random measures}\label{sec:pt}
An important class of random counting measures are the Poisson-type (PT; Poisson, binomial, and negative binomial). They are uniquely characterized by distributional closure under restriction (binomial thinning). For a given PT random measure, the counts in all subspaces (volume elements) are PT. It turns out that the PT distributions are the only members of the power series family of distributions to possess this self-similarity property. For completeness we give the result, which also illustrates restrictions (traces).

\begin{samepage}
\begin{thm}[Existence and Uniqueness of PT Random Measures]\label{thm:ptrms} Assume that $K\sim \kappa_\theta$ where pgf $\psi_\theta$ belongs to the canonical non-negative power series family of distributions and $\{0,1\}\subset supp(K)$. Consider the random measure $N=(\kappa_\theta,\nu)$ on the  space $(E, \mathscr{E})$ and assume that $\nu$ is diffuse. Then for any $A\subset E$ with $\nu(A)=a>0$ there exists a mapping $h_a:\Theta\rightarrow \Theta$ such that the restricted random measure is $N_A=(\kappa_{h_a(\theta)}, \nu_A)$, that is, \begin{equation*}\E e^{-N_Af}=\psi_{h_a(\theta)}(\nu_A e^{-f})\quad \text{for}\quad f\in \mathscr{E}_{\ge0}\end{equation*}  iff $K$ is Poisson, negative binomial or binomial. 
\end{thm}
\end{samepage}\begin{proof} Theorem 3 \cite{Bastian:2020aa}.\end{proof}

Note that if $\nu$ is atomic, then the `if' holds, so that the PT distributions in subsets are rescaled versions of themselves, but the `only if' (uniqueness) does not. Thus the PT random measures are closed under restriction to all subspaces, uniquely so in the power series family of distributions for diffuse measures.

\subsection{Random fields from random measures}\label{sec:rf}

A non-negative random field $G=\{G(y): y\in F\}$ on $(F,\mathscr{F})$ is a collection of random variables $G(y)$ such that the mapping $(\omega,y)\mapsto G(\omega,y)$ is measurable relative to $\mathscr{H}\otimes\mathscr{F}$ and $\mathscr{B}_{\R_{\ge0}}$. In view of Definition \ref{de:rt}, $G$ may be interpreted as a non-negative random transformation. We can construct non-negative random fields from random measures as follows. Let $N$ be a random measure on $(E,\mathscr{E})$, e.g., STC or FARW. Let $k:E\times F\mapsto\R_{\ge0}$ be a non-negative $\mathscr{E}\otimes\mathscr{F}$-measurable function. Then \begin{equation}\label{eq:rf}G(y) = \int_E N(\D x)k(x,y)\for y\in F\end{equation} forms a non-negative random field $G$ on $(F,\mathscr{F})$. The properties of $G$ follow from $N$, e.g., $\E G(y) = \E Nk(\cdot,y)$.

\section{STC random measure graph models}\label{sec:rg}Here we apply the results of the previous section to graph (network) models. The idea is to use a product random measure to represent the edge set and the composition of the test function and transformation to define a non-negative edge weight function. The edge weight function is considered to be deterministic or random. This representation is prototypical to labeled random graphs.
  
\begin{samepage}
\begin{de}[STC random graph]\label{def:rg} Consider the STC random measure $N=(\kappa,\nu)$ on $(E,\mathscr{E})$ formed by independency $\mathbf{X}=\{X_i: i=1,\dotsb,K\}$. Now consider the product random measure $M=N\times N$ on $(E\times E,\mathscr{E}\otimes\mathscr{E})$ formed by $\mathbf{X}\times\mathbf{X}=\{(X_i,X_j): i=1,\dotsb,K,j=1,\dotsb,K\}$. For (possibly random) transformation $f:E\times E\mapsto F$, where $(F,\mathscr{F})$ is a measurable space, and function $g\in\mathscr{F}_{\ge0}$, the triple $(M,f,g)$ defines a non-negative weighted labeled \emph{random graph} \[G(g\circ f)=(\mathbf{X},\mathbf{X}\times\mathbf{X},g\circ f)\] where $\mathbf{X}$ is the labeled vertex set, $\mathbf{X}\times\mathbf{X}$ is the labeled edge set, and $g\circ f\in(\mathscr{E}\otimes\mathscr{E})_{\ge0}$ is the labeled edge weight function. The underlying vertex and edge sets are $\mathbf{V}=\{1,\dotsb,K\}$ and $\mathbf{E}=\{1,\dotsb,K\}^2$ respectively. \end{de}
\end{samepage}

\begin{re}[Notation]\label{re:notate} For $(F,\mathscr{F})=(\R_{\ge0},\mathscr{B}_{\R_{\ge0}})$ and identity function $I(x)=x$, we have $G(f)=G(I(f))$. 
\end{re}

Next we give a series of definitions for the support, edge number / weight, vertex number, isolated components, and active graph of $G(g\circ f)$. 

\begin{de}[Support]The support of $G=G(g\circ f)$ is the labeled edge set \[\supp(G)\equiv\{g\circ f>0\}=\{(x,y)\in\mathbf{X}\times\mathbf{X}: g\circ f(x,y)>0\}\] We refer to such edges as active or non-trivial. \end{de}

\begin{de}[Active edge number/weight]The number of active edges of $G=G(g\circ f)$ is $e(G)\equiv|\supp(G)|=(M\circ f^{-1})\ind{\supp(g)}$ and the active edge weight of $G$ is $w(G)\equiv(M\circ f^{-1})g$. Their means are $\E e(G) = \E(M\circ f^{-1})\ind{\supp(g)}$ and $\E w(G) = \E (M\circ f^{-1})g$ respectively. If $g\circ f$ is symmetric, then $G(g\circ f)$ is undirected, and graph active edge count and weight are inflated by the double counting of $M$. In this case, the edge count and weight for symmetric $g\circ f$ are defined in terms of the normalized $M$, $M^*$ (Remark~\ref{re:norm}), as $e^*(G)\equiv(M^*\circ f^{-1})\ind{\supp(g)}$ and $w^*(G)\equiv(M^*\circ f^{-1})g$ with respective means $\E e^*(G) = \E(M^*\circ f^{-1})\ind{\supp(g)}$ and $\E w^*(G) = \E (M^*\circ f^{-1})g$.
\end{de}

\begin{de}[Dense/sparse] Assuming $c$ is of greater order than $\delta$, the random graph model $G=G(g\circ f)$ is said to be \emph{dense} if $\E e(G)$ grows as $O(c^2)$ and is said to be \emph{sparse} if $\E e(G)$ grows as $o(c^2)$.  
\end{de}

\begin{de}[Active vertices] The set of active (non-trivial) vertices of $G=G(g\circ f)$ is \[t(G) = \{x\in\mathbf{X}: \text{for any }y\in\mathbf{X},\, (x,y)\text{ or }(y,x)\in\supp(G)\}\subseteq\mathbf{X}\] and the number of active vertices is $v(G) = |t(G)|$.
\end{de}

\begin{de}[Isolated components] The isolated components of $G(g\circ f)$ correspond to its inactive (trivial) vertices.
\end{de}
 
\begin{de}[Active graph] The active graph $H$ of random $G=G(g\circ f)$ is the random subgraph of $G$ formed as $H(g\circ f)=(t(G),t(G)\times t(G),g\circ f)$. 
\end{de}

\begin{re}[Active graph] By construction the active graph $H$ of $G$ and $G$ have identical active edge number and weight: $e(H)=e(G)$ and $w(H)=w(G)$.
\end{re}

The definitions of $f:E\times E\mapsto F$ and $g:F\mapsto\R_{\ge0}$ evidently encode the properties of $G(g\circ f)$. In most of what follows we assume $(F,\mathscr{F})=(\R_{\ge0},\mathscr{B}_{\R_{\ge0}})$, identity $g(x)=I(x)=x$, and deterministic $f$, which gives the random graph $G(f)$. When $f(x,x)=0$ for $x\in E$, self-edges are avoided, and we point this out throughout the following whether this is in-force or not. If $f$ is symmetric, then $G(f)$ is undirected. If $f$ is symmetric and $\{0,1\}$-valued with $f(x,x)=0$, then $G(f)$ is simple and undirected. And so and so forth. The number of underlying vertices is $K$ with mean $c$, and the number of underlying edges is $K^2$ with mean $c^2+\delta^2$. Note that symmetric $f$ taking values in $[0,1]$ is known in graph theory as a \emph{graphon} (graph-function) \citep{graphon}, typically with $E=[0,1]$, $\nu=\Leb$, and $f(x,x)=0$ for $x\in E$ (to avoid self-edges). For atomic $\nu$, the labeled edge collection $\mathbf{X}=\{X_i:i=1,\dotsb,K\}$ is a \emph{multiset} (or disjoint union) of $K$ iid discrete random variables and hence the edge collection is also a multiset. 



In view of Theorem~\ref{thm:ptrms}, when $N=(\kappa,\nu)$ is Poisson-type (PT), i.e. $\kappa_\theta$ is Poisson, negative binomial, or binomial, with canonical parameter $\theta\in\Theta$ and pgf $\psi_\theta$, then the pgf of $K_A=N\ind{A}$ for $A\subseteq E$ with $a=\nu(A)>0$ is given by $\psi_A(t)=\psi_{h_a(\theta)}(t)$, where $h_a:\Theta\mapsto\Theta$ is a rescaling mapping. More generally, for restricted measure $N_A=(\kappa_{h_a(\theta)},\nu_A)$, where $\nu_A(\cdot)=\nu(A\cap\cdot)/a$ is the restriction of $\nu$ to $A$, and test function $f\in\mathscr{E}_{\ge0}$, the Laplace transform of $N_Af$ is \[\varphi(\alpha)=\psi_{h_a(\theta)}(\nu_A e^{-\alpha f})\for\alpha\in\R_{\ge0}\] Hence properties derived from the PT random measure $N$ acting on test functions or under transformation, such as degree, inherit this distributional self-similarity property, i.e. closure under restriction to subspaces.



This remainder of this section is organized as follows. Next we discuss random graphs $G(f)=G(I\circ f)$ for deterministic functions $f\ge0$ (\ref{sec:deter}). We follow this for random graphs $G(g\circ f)$ having random transformations $f$ (\ref{sec:rand}) and remark on their generalized spectral (\ref{sec:svd}) and Sobol representations (\ref{sec:sobol}). We then discuss triangle functions (\ref{sec:triangle}), degree functions (\ref{sec:degree}), and degree distributions (\ref{sec:degreedis}) for deterministic and random transformations. Next we give a criterion for existence of giant components (\ref{sec:giant}) and treat active vertices (\ref{sec:part}). Finally, we give some particular kinds of graph: empirical (\ref{sec:emp}), directed acyclic (\ref{sec:dag}), and rewired (\ref{sec:resample}). 

\subsection{Deterministic transformations}\label{sec:deter} Here we consider the random graph $G(f)=(\mathbf{X},\mathbf{X}\times\mathbf{X},f)$ for deterministic $f\in(\mathscr{E}\otimes\mathscr{E})_{\ge0}$. By Theorem~\ref{thm:product}, the product random measure $M$ has mean \[\E Mf = c(\nu\times I)f + (c^2+\delta^2-c)(\nu\times\nu)f\] For the indicator function $f=\ind{A}$ of subspace $A\subseteq E\times E$, then the \emph{expected number of edges} of $G(\ind{A})$ is given by \[\E e(G) = \E M(A) = c(\nu\times I)(A) + (c^2+\delta^2-c)(\nu\times\nu)(A)\] and we have $\E M(E\times E)=\E K^2 = c^2+\delta^2$. The first term $\E e_{\text{self}}(G)=c(\nu\times I)(A)$ is the mean number of \emph{self-edges} of $G(\ind{A})$, whereas the second term $\E e_{\text{ext}}(G)=(c^2+\delta^2-c)(\nu\times\nu)(A)$ is the mean number \emph{external edges} of $G(\ind{A})$, with $\E e(G) = \E e_{\text{self}}(G)+ \E e_{\text{ext}}(G)$. 

The indicator function $\ind{A}$ encodes the restricted product random measure $M_A=M\ind{A}$ with mean measure $\E M_A = \E M\ind{A}$. We think of the set $A$ as the \emph{interaction space} of the graph $G(\ind{A})$. Put $a=(\nu\times\nu)(A)$. Assuming $c^2$  on order dominates $\delta^2$, then the mean number of edges of $G(\ind{A})$ scales $O(ac^2)$. If $a$ is fixed, then the scaling is $O(c^2)$ and the graph is said to be \emph{dense}. When the number of edges is $o(c^2)$, the graph is \emph{sparse}. If $a\sim1/c$, then the scaling is $O(c)$ and the graph $G(\ind{A})$ is sparse. Thus the choice of $N=(\kappa,\nu)$ and $f=\ind{A}$ determines the graph regime.

To avoid self-edges, we choose $A$ through the following recipe: take arbitrary $B\subseteq E\times E$, $C=\{(x,y)\in E\times E: x\ne y\}$, and set $A=B\cap C$; then, $\E M(A)=(c^2+\delta^2-c)(\nu\times\nu)(A)$ is the mean number of edges, excluding self-edges.


\subsection{Random transformations}\label{sec:rand}

Here the edge weight function is a random transformation with a given marginal transition kernel. We describe four random transformations: Bernoulli, digraphon, binomial, and Poisson. They are $\{0,1\}$, $\{0,1\}^2$, $\{0,1,\dotsb,n\}$, and $\N_{\ge0}$ valued respectively. There are many more. All the results follow from Theorem~\ref{thm:product} on the mean measure of the product STC random measure.


\subsubsection{Bernoulli transformation} Consider the random graph $G(\phi)=(\mathbf{X},\mathbf{X}\times\mathbf{X},\phi)$ formed from $(M,\phi,I)$, with random transformation $\phi:E\times E\mapsto \R_{\ge0}$ independent of $M$ with marginal transition kernel $Q$. The image random measure of $M$ under $\phi$ on $(F,\mathscr{F})=(\R_{\ge0},\mathscr{B}_{\R_{\ge0}})$ is given by \[(M\circ\phi^{-1})f = \sum_i^K\sum_j^Kf(\phi(X_i,X_j))\for f\in\mathscr{F}_{\ge0}\] By Theorem~\ref{thm:product}, $M\circ\phi^{-1}$ has mean \[\E(M\circ\phi^{-1})f =c(\nu\times I)Qf + (c^2+\delta^2-c)(\nu\times\nu)Qf\for f\in\mathscr{F}_{\ge0}\] Note that $\phi(x,x)=0$ is encoded by $Q((x,x),\cdot)=\delta_0(\cdot)$. The expected edge weight of $G(f(\phi))$ is given by $\E(M\circ\phi^{-1})f$. When $\phi$ is $\{0,1\}$-valued, then the expected edge weight is the expected number of edges. We have the following theorem for Bernoulli $\phi$. Note that neither the transformation $\phi$ nor the function $f$ is necessarily symmetric.

\begin{samepage}
\begin{thm}[Bernoulli] Let $M=N\times N$ be a product random measure on $(E\times E,\mathscr{E}\otimes\mathscr{E})$ where $N=(\kappa,\nu)$. Consider a random transformation $\phi:E\times E\mapsto\{0,1\}$ independent of $M$ with marginal transition kernel $Q((x,y),\cdot)=\text{Bernoulli}(f(x,y))$ for $f\in(\mathscr{E}\otimes\mathscr{E})_{[0,1]}$. Then \[\E(M\circ\phi^{-1})I = \E Mf = c(\nu\times I)f + (c^2+\delta^2-c)(\nu\times\nu)f\]  where $I$ is the identity function $I(z)=z$ or kernel $I(z,\cdot)=\delta_z(\cdot)$. 
\end{thm}
\end{samepage}
\begin{proof}\begin{align*}\E(M\circ\phi^{-1})I &=c(\nu\times I)QI + (c^2+\delta^2-c)(\nu\times\nu)QI\\&= c\int_{E\times \{0,1\}}\nu(\D x)Q((x,x),\D z)z + (c^2+\delta^2-c)\int_{E\times E\times\{0,1\}}\nu(\D x)\nu(\D y)Q((x,y),\D z)z\\&=c\int_{E}\nu(\D x)f(x,x) + (c^2+\delta^2-c)\int_{E\times E}\nu(\D x)\nu(\D y)f(x,y)\\&=c(\nu\times I)f + (c^2+\delta^2-c)(\nu\times\nu)f\\&=\E Mf\end{align*}
\end{proof}

\begin{re}[Bernoulli edge count] For Bernoulli random graph $G=G(\phi)$ with $f\in(\mathscr{E}\otimes\mathscr{E})_{[0,1]}$, the mean number of edges is $\E e(G) = \E Mf$. 
\end{re}

\begin{re}[Graph regime]For Bernoulli random transformations $\phi$ built from $f\in (\mathscr{E}\otimes\mathscr{E})_{[0,1]}$, the structure of $f$ determine the graph regime of sparse or dense, where $a=(\nu\times\nu)f$.
\end{re}
 
\begin{re}[Weighted Bernoulli graphs]To construct weighted graphs based on the Bernoulli transformation $\phi:E\times E\mapsto\{0,1\}$ and $f\in(\mathscr{E}\otimes\mathscr{E})_{[0,1]}$---weighted Bernoulli graphs---we introduce a deterministic positive weight function $g\in(\mathscr{E}\otimes\mathscr{E})_{>0}$ to define a non-negative random transformation $\phi_g:E\times E\mapsto\R_{\ge0}$ as $\phi_g=g\phi$. The $(M,\phi_g,I)$ forms the weighted Bernoulli graph $G(\phi_g)$ with mean edge count $\E(M\circ\phi_g^{-1})\ind{\supp(I)}=\E Mf$ and weight $\E(M\circ\phi_g^{-1})I = \E M(fg)$.
\end{re}

\paragraph{Relation to the graphex} The class of random graphs as formulated $(M,f,g)$ in Definition~\ref{def:rg} is similar to, for example in the undirected case for certain $f$ and $g$, a \emph{graphex} \citep{graphex}. A graphex parameterizes a family of random graphs whose vertices are in $\R_{\ge0}$ and whose active edge set is defined through a finite symmetric jointly exchangeable point process $\Gamma$ on $\R_{\ge0}^2$, constructed from a latent Poisson process $\Pi=\{(\theta_i,\vartheta_i): i\ge1\}$ on $\R_{\ge0}^2$. Formally a graphex is a triple $(I,S,W)$ on a measure space $(E,\mathscr{E},\nu)$ where $I\in\R_{\ge0}$ encodes the (Poisson) rate of isolated edges, integrable $S:E\mapsto\R_{\ge0}$ encodes the (Poisson) dependent rate of star components, and symmetric $W:E\times E\mapsto[0,1]$ is a graphon. The graphex $(0,0,W)$ construction of $\Gamma$ is as follows: for every distinct pair of points $(\theta,\vartheta),(\theta',\vartheta')\in\Pi$, the edge $(\theta,\theta')$ belongs to $\Gamma$ with probability $W(\vartheta,\vartheta')$. Then the class of graphs generated by the graphex $(0,0,W)$ is equivalent to the active graphs of $(M,f,I)$ in the case of Poisson $\kappa$ and symmetric random transformation $f:E\times E\mapsto\{0,1\}$ having (marginal) transition kernel $Q((x,y),\cdot)\sim\text{Bernoulli}(W(x,y))$. 

\subsubsection{Bernoulli examples} We give a some examples of Bernoulli random transformations $\phi$. 

\paragraph{Example: Erd\"{o}s-Renyi (ER)} Let $\phi$ be symmetric and $\{0,1\}$-valued with $\phi(x,x)=0$ following from $f(x,x)=0$. The ER graph $G(n,p)$ may be attained by taking $\kappa=\text{Dirac}(n)$ and defining the marginal transition kernel as $Q((x,y),\cdot)=\text{Bernoulli}(p)$ for $x\ne y$ and noting that $f(x,x)=0$ encodes $Q((x,x),\cdot)=\delta_0(\cdot)$. Then \[G(n,p) = G(\phi)\] The symmetry and support of $\phi$ imply that there are $\frac{1}{2}n(n-1)=\binom{n}{2}$ independent $
\text{Bernoulli}(p)$ random variables that define a random realization of $\phi$. The mean number of active edges is $\E e(G)=\frac{1}{2}n(n-1)p = \binom{n}{2}p$.



\paragraph{Example: Power-law} Let $E=[0,1]$ and $\nu=\Leb$. Define power-law graphon $f\in(\mathscr{E}\otimes\mathscr{E})_{[0,1]}$  as $f(x,y)=\ind{}(x\ne y)(1+bx)^{-2}(1+by)^{-2}$ for $b>0$. Then \[\E(M\circ\phi^{-1})I = (c^2+\delta^2-c)(\nu\times\nu)f = (c^2+\delta^2-c)/(1+b)^2\] For $b=\sqrt{c}$ and assuming $c^2$ is of greater order than $\delta^2$, then $\E(M\circ\phi^{-1})I\sim c$ as $c\rightarrow\infty$. This parameterization yields a sparse graph model.

\paragraph{Example: Exponential} Let $E=[0,1]$ and $\nu=\Leb$. Define exponential graphon $f\in(\mathscr{E}\otimes\mathscr{E})_{[0,1]}$  as $f(x,y)=\ind{}(x\ne y)\exp_-b(x+y)$ for $b>0$. Then \[\E(M\circ\phi^{-1})I = (c^2+\delta^2-c)(\nu\times\nu)f = (c^2+\delta^2-c)\left(\frac{1-e^{-b}}{b}\right)^2\] As with power-law, for $b=\sqrt{c}$ and assuming $c^2$ is of greater order than $\delta^2$, then $\E(M\circ\phi^{-1})I\sim c$ as $c\rightarrow\infty$ yields a sparse graph model.

\paragraph{Example: Block} Let $A_1,\dotsb,A_n$ be a partition of $E$. For every $(i,j)\in\{1,\dotsb,n\}^2$, let $p_{ij}\ge0$ be a probability. Define  block graphon $f(x,y)=\sum_i\sum_j p_{ij}\ind{A_i\times A_j}(x,y)$. Then \[\E(M\circ\phi^{-1})I = c\sum_ip_{ii}\nu(A_i) + (c^2+\delta^2-c)\sum_i\sum_jp_{ij}(\nu\times\nu)(A_i\times A_j)\] This is a \emph{stochastic block model} of random graphs, where the partition is the collection of communities of the vertices \citep{block}.

\paragraph{Example: Dot-product} Let $E=[0,1]^d$ and take $\nu=\Leb$. Define dot-product graphon $f\in(\mathscr{E}\otimes\mathscr{E})_{[0,1]}$ as $f(x,y) = \langle x^a,y^a\rangle / d$ where $x^a\equiv(x_1^a,\dotsb,x_d^a)$ for $a\ge0$. Then \[\E(M\circ\phi^{-1})I = \frac{c^2+\delta^2+ac}{(a+1)^2}\] is free of $d$. This is the \emph{random dot-product model} of random graphs \citep{dot}.



\subsubsection{Digraphon transformation} Here we consider the triple $(M,\phi,f)$ and its random graph $G(f\circ\phi)$ for $(F,\mathscr{F})=(\{0,1\}^2,2^F)$. The random transformation $\phi$ is based on the \emph{digraphon} $f=(f_{00},f_{01},f_{10},f_{11},g)$, where \begin{enumerate}\item[(i)]$f_{00},f_{01},f_{10},f_{11}\in(\mathscr{E}\otimes\mathscr{E})_{[0,1]}:\quad f_{00}+f_{01}+f_{10}+f_{11}=1$\item[(ii)] $f_{00}$ and $f_{11}$ are symmetric \item[(iii)] $f_{01}(x,y)=f_{10}(y,x)$ for $(x,y)\in E\times E$\item[(iv)] $g\in\mathscr{E}_{[0,1]}$\end{enumerate} \citep{digraphon}. For each digraphon $(f_{00},f_{01},f_{10},f_{11},g)$, the $(f_{00},f_{01},f_{10},f_{11})$ defines a transition kernel $Q$ from $E\times E$ into $F$ as \[Q((x,y),(a,b)) = f_{ab}(x,y)\] so that \[Q((x,y),\cdot)\sim\sum_{(a,b)\in F}f_{ab}(x,y)\] We have the following lemma on the symmetry of $Q$. 

\begin{lem}[Symmetry]\label{lem:sym} $Q((x,y),A)=Q((y,x),A)$ for $(x,y)\in E\times E$ and $A\subseteq F$ such that either $(0,1),(1,0)\in A$ or $(0,1),(1,0)\notin A$.
\end{lem}
\begin{proof} The result follows from the symmetries of $f_{00}, f_{01}+f_{10},f_{11}$ and their sums. 
\end{proof}

Let $R:E\times E\mapsto F$ be a transition kernel \[R((x,x),\cdot)\sim\text{Bernoulli}(g(x))\times I\] and define the transition kernel $P:E\times E\mapsto F$ as \begin{equation}\label{eq:P} P((x,y),\cdot)=\begin{cases} R((x,x),\cdot) &\text{if }x=y\\Q((x,y),\cdot)&\text{otherwise}\end{cases}\end{equation} such that $(\phi(x,y),\phi(y,x))\sim P((x,y),\cdot)$. Then the mean of $M\circ\phi^{-1}$ is \begin{align*}\E(M\circ\phi^{-1})f &= (\E M)Pf \\&= c(\nu\times I)Rf + (c^2+\delta^2-c)(\nu\times\nu)Qf\\&=c\sum_{(a,b)\in F}\int_E\nu(\D x)R((x,x),(a,b))f(a,b) + (c^2+\delta^2-c)\sum_{(a,b)\in F}\int_{E\times E}\nu(\D x)\nu(\D y)Q((x,y),(a,b))f(a,b)\\&=c\nu(\D x)((1-g(x))f(0,0)+g(x)f(1,1))+(c^2+\delta^2-c)\sum_{(a,b)\in F}\int_{E\times E}\nu(\D x)\nu(\D y)f_{ab}(x,y)f(a,b)\end{align*} Organizing these, we have the following result. 

\begin{samepage}
\begin{thm}[Digraphon] Let $M=N\times N$ be a product random measure on $(E\times E,\mathscr{E}\otimes\mathscr{E})$ where $N=(\kappa,\nu)$. Put $(F,\mathscr{F})=(\{0,1\}^2,2^F)$. Consider a random transformation $\phi:E\times E\mapsto F$ independent of $M$ with marginal transition kernel $P$ \eqref{eq:P}. Then \[\E(M\circ\phi^{-1})f = c\nu((1-g)f(0,0)+gf(1,1)) + (c^2+\delta^2-c)\sum_{(a,b)\in F}(\nu\times\nu)f_{ab}f(a,b)\for f\in\mathscr{F}_{\ge0}\] 
\end{thm}
\end{samepage}

Consider $C=\{(0,1),(1,0),(1,1)\}$. Then $\E e_C = \E(M\circ\phi^{-1})(C) = c\nu(g)+(c^2+\delta^2-c)(\nu\times\nu)(f_{01}+f_{10}+f_{11})$ is the mean number of edges of $G(\ind{C}(\phi))$. Defining $A=\{(0,1),(1,0)\}$, the mean number of directed edges of $G(\ind{A}(\phi))$ is given by $\E e_A=\E(M\circ\phi^{-1})(A)=(c^2+\delta^2-c)(\nu\times\nu)(f_{01}+f_{10})$. Defining $B=\{(1,1)\}$, the mean number of bidirected edges of $G(\ind{B}(\phi))$ is given by $\E e_B = \E(M\circ\phi^{-1})(B)=c\nu(g)+(c^2+\delta^2-c)(\nu\times\nu)f_{11}$. 

\subsubsection{Binomial transformation} Consider the random graph $G(\phi)=(\mathbf{X},\mathbf{X}\times\mathbf{X},\phi)$ formed from $(M,\phi,I)$ with binomial transformation $\phi$ having marginal transition kernel $Q((x,y),\cdot)=\text{Binomial}(n,w(x,y))$ for $n\ge1$ and $w\in(\mathscr{E}\otimes\mathscr{E})_{[0,1]}$. This transformation is $\{0,1,\dotsb,n\}$-valued and results in a multigraph for $n>1$. Of course, when taking $n=1$, the Bernoulli transformation is retrieved.

\begin{samepage}
\begin{thm}[Binomial] Let $M=N\times N$ be a product random measure on $(E\times E,\mathscr{E}\otimes\mathscr{E})$ where $N=(\kappa,\nu)$. Consider a random transformation $\phi:E\times E\mapsto\{0,1,\dotsb,n\}$ independent of $M$ with marginal transition kernel $Q((x,y),\cdot)=\text{Binomial}(n,f(x,y))$ for $f\in(\mathscr{E}\otimes\mathscr{E})_{[0,1]}$. Then \[\E(M\circ\phi^{-1})I = \E Mf = nc(\nu\times I)f + n(c^2+\delta^2-c)(\nu\times\nu)f\]  where $I$ is the identity function $I(z)=z$ or kernel $I(z,\cdot)=\delta_z(\cdot)$. 
\end{thm}
\end{samepage}

\subsubsection{Poisson transformation} Consider the random graph $G(\phi)=(\mathbf{X},\mathbf{X}\times\mathbf{X},\phi)$ formed from $(M,\phi,I)$ with Poisson random transformation $\phi$ having marginal transition kernel $Q((x,y),\cdot)=\text{Poisson}(f(x,y))$ for $f\in(\mathscr{E}\otimes\mathscr{E})_{\ge0}$. This transformation is non-negative integer-valued and results in a multigraph. Note that for small $f(x,y)$, the Poisson transformation is approximately Bernoulli, $\phi(x,y)\simeq\text{Bernoulli}(f(x,y))$.

\begin{samepage}
\begin{thm}[Poisson] Let $M=N\times N$ be a product random measure on $(E\times E,\mathscr{E}\otimes\mathscr{E})$ where $N=(\kappa,\nu)$. Consider a random transformation $\phi:E\times E\mapsto\N_{\ge0}$ independent of $M$ with marginal transition kernel $Q((x,y),\cdot)=\text{Poisson}(f(x,y))$ for $f\in(\mathscr{E}\otimes\mathscr{E})_{\ge0}$. Then \[\E(M\circ\phi^{-1})I = \E Mf = c(\nu\times I)f + (c^2+\delta^2-c)(\nu\times\nu)f\]  where $I$ is the identity function $I(z)=z$ or kernel $I(z,\cdot)=\delta_z(\cdot)$. 
\end{thm}
\end{samepage}

\subsection{Generalized spectral representation}\label{sec:svd} Consider the graph $G(w=g\circ f)$ for random $f$ (and $w$). Denote mean weight $W = \E w$. The mean weight $W\in L^2(\nu\times\nu)$ forms a compact operator $T_{W}$ on the Hilbert space $L^2(\nu)$ (c.f. Theorem~\ref{re:spectral}) with generalized spectral representation (singular value decomposition) \[W(x,y) = \sum_n\sigma_nf_n(x)g_n(y) \for (x,y)\in E\times E\] Such theory is well developed for and commonly applied to random graphs in the case of Bernoulli transformations $w$ based on graphons $W$, i.e. symmetric $W\in(\mathscr{E}\otimes\mathscr{E})_{[0,1]}$, where $T_{W}$ is self-adjoint \citep{graphon}. When $W$ is symmetric, another approach is to use the normalized graph Laplacian kernel $p(x,y) = W(x,y) / d(x)$, where $d(x)=(T_{W}\ind{E})(x) = \int_E W(x,y)\nu(\D y)$, which for finite graphs yields the diffusion map \citep{diffusionmap}. 

\subsection{Sobol system representation}\label{sec:sobol} 

A classical \emph{Sobol system} \citep{sobol93}, also known as a high dimensional model representation \citep{alis99} or a functional ANOVA expansion \citep{hooker}, is a representation of a square-integrable function of $n$-variates on a product space as a superposition of $2^n$ orthogonal component functions (projections). In particular, consider the Hilbert space $\mathcal{H}=L^2(E\times E,\mathscr{E}\otimes\mathscr{E},\nu\times\nu)$ with inner product $\langle\cdot,\cdot\rangle_{\mathcal{H}}$. The Sobol system representation of $\mathcal{H}$ is a direct sum of orthogonal subspaces \[\mathcal{H} = \mathcal{V}_0\oplus\mathcal{V}_1\oplus\mathcal{V}_2\oplus\mathcal{V}_{12}\] defined as \begin{align*}\mathcal{V}_0&=\{f\in\mathcal{H}: f=C \text{ for constant }C\in\R\}\\\mathcal{V}_1 &= \{f\in \mathcal{H}: f=f_1(x)\text{ with }\nu f_1=0\}\\\mathcal{V}_2 &= \{f\in \mathcal{H}: f=f_2(y)\text{ with }\nu f_2=0\}\\\mathcal{V}_{12}&=\{f\in\mathcal{H}: f=f_{12}(x,y)\text{ with }\nu f_{12}(x,\cdot)=\nu f_{12}(\cdot, y)=0\}\end{align*}


Let $\mathcal{W}=(\mathscr{E}\otimes\mathscr{E})_{\ge0}\cap\mathcal{H}$ be the subcollection of non-negative square-integrable functions. Consider the random graph $G(w=g\circ f)$ with mean edge weight function $W=\E w \in\mathcal{W}$. Its Sobol system is given by \[W(x,y) = W_0 + W_1(x)+W_2(y)+W_{12}(x,y)\for (x,y)\in E\times E\] where $(W_0, W_1, W_2, W_{12})\in\mathcal{V}_0\otimes\mathcal{V}_1\otimes\mathcal{V}_2\otimes\mathcal{V}_{12}$ are the component functions (orthogonal projections) of $W$ constructed as \begin{align*}W_0&\equiv(\nu\times\nu)W\\W_1(x) &\equiv\nu W(x,\cdot) - W_0\\W_2(y)&\equiv\nu W(\cdot,y) - W_0\\W_{12}(x,y)&\equiv W(x,y)-W_0-W_1(x)-W_2(y)\end{align*} The mutually orthogonal $W_1,W_2,W_{12}$ convey a decomposition of variance (ANOVA) \[\Var W = \Var W_1 + \Var W_2 + \Var W_{12}\] The ANOVA reveals the independent and cooperative coordinate contributions to the variance. When $\Var W>0$, the normalized variances $\amsbb{S}_u\equiv\Var W_u / \Var W\in[0,1]$ are referred to as \emph{sensitivity indices}. 

The \emph{effective dimension} of $W$ is a measure of the order complexity of $W$ and is defined as \[ED(W)\equiv \begin{cases}1\times(\amsbb{S}_1+\amsbb{S}_2)+2\times\amsbb{S}_{12}&\text{if }\Var W >0\\0&\text{otherwise}\end{cases}\] 

For example, consider $E=[0,1]$, $\nu=\Leb$, and $W(x,y) = \exp_-a(x+y)$. The sensitivity indices \begin{align*}\amsbb{S}_1=\amsbb{S}_2&=\frac{2(e^a - 1)}{(a+2)e^a+a-2}\xrightarrow{a\rightarrow\infty}0\\\amsbb{S}_{12}&=\frac{(a-2)e^a + a + 2}{(a+2)e^a+a-2}\xrightarrow{a\rightarrow\infty}1\end{align*} reveal that the effective dimension of $W$ goes to two as $a$ goes to infinity.

The normalized mean in- and out-degree functions are retrieved from the integral operators used to define the first-order component functions  as \begin{align}D_o(x) &\equiv \nu W(x,\cdot)=W_0+W_1(x)\label{eq:Do}\\D_i(x) &\equiv \nu W(\cdot,x)= W_0+W_2(x)\label{eq:Di}\end{align} This gives \begin{align*}W(x,y) &=W_0+W_1(x)+W_2(y)+W_{12}(x,y)\\&=\widetilde{W}(x,y) + W_{12}(x,y)\\&= D_o(x)+ D_i(y) - W_0+W_{12}(x,y)\end{align*} where $\widetilde{W}$ is the first-order Sobol representation of $W$. $\widetilde{W}$ encodes the total weight and normalized degree functions \[W_0=(\nu\times\nu)\widetilde{W},\quad D_o(x)=\nu\widetilde{W}(x,\cdot), \quad D_i(x)=\nu\widetilde{W}(\cdot,x)\] and can take negative values. 


The second-order term $W_{12}$ captures the contribution of coordinate interactions to $W$. It is orthogonal to $W_0$, $W_1$, and $W_2$ and their superpositions and hence to $\widetilde{W}$ and the degree functions $D_i$ and $D_o$. This suggests the notion of additive weight functions.

\begin{de}[Additivity]We refer to a weight function $W\in\mathcal{W}$ as \emph{additive} when $W_{12}=0$.\end{de} 
\begin{prop}[Additivity] $W\in\mathcal{W}$ is additive iff $W(x,y)=c+W_a(x)+W_b(y)$ for some $c\ge0$ and $W_a,W_b\in L^2(E,\mathscr{E},\nu)$.\end{prop}

Consider the additive weight function $W(x,y)=W_a(x)+W_b(y)$ for $W_a,W_b\in\mathscr{E}_{\ge0}\cap L^2(E,\mathscr{E},\nu)$. It has component functions \[W_0= \nu(W_a)+\nu(W_b),\quad W_1(x) = W_a(x)-\nu(W_a),\quad W_2(y) = W_b(y)-\nu(W_b)\] and degrees  \[D_o(x) =W_a(x)+\nu(W_b),\quad D_i(y) =W_b(y)+\nu(W_a)\]





\subsection{Triangle function}\label{sec:triangle} Consider a random graph $G(g\circ f)$ formed from $(M,f,g)$ with random transformation $f$. We define its triangle (weight) function as the sum of the product of weights of distinct vertex triplets.

\begin{de}[Triangle function] The triangle function of the random graph $G(g\circ f)$ is defined as \[\triangle(z) = \int_{E\times E}M(\D x,\D y)\ind{}(x\ne y\ne z)g\circ f(x,y)\, g\circ f(y,z)\, g\circ f(z,x)\for z\in E\]
\end{de} 

\begin{re}[Random field] The triangle function forms a random field, i.e. it embodies \eqref{eq:rf}. Putting \[k_z(x,y)\equiv k(x,y,z)\equiv\ind{}(x\ne y\ne z)g\circ f(x,y)\, g\circ f(y,z)\, g\circ f(z,x)\for (x,y,z)\in E\times E\times E\] then $\triangle(z)=Mk_z$.
\end{re}

\begin{prop}[Mean] For $W=\E(g\circ f)$ and $\widetilde{W}(x,y)\equiv\ind{}(x\ne y)W(x,y)$, then \[\E\triangle(z) = (c^2+\delta^2-c)\int_{E\times E}\nu(\D x)\nu(\D y)\widetilde{W}(x,y)\widetilde{W}(y,z)\widetilde{W}(z,x)\for z\in E\]
\end{prop}

\subsection{Degree function}\label{sec:degree} An important object is the degree function $d:E\mapsto\R_{\ge0}$ of a random graph $G(g\circ f)$. We give its definition below and an immediate remark.

\begin{de}[Degree functions] The out-degree and in-degree functions of the random graph $G(g\circ f)$ are respectively defined as \[d_o(x) \equiv \int_EN(\D y)g\circ f(x,y)\for x\in E\] and \[d_i(x) \equiv \int_EN(\D y)g\circ f(y,x)\for x\in E\] Writing $f_x(y)=f(x,y)$ and $f^x(y)=f(y,x)$, these become $d_o(x) = (N\circ f_x^{-1})g$ and $d_i(x)=(N\circ f^{x-1})g$.
\end{de}

\begin{re}[Random field] The degree functions form random fields.
\end{re}

\begin{prop}[Means] Let $w=g\circ f$ have mean $W=\E w$. Then for $x\in E$ \begin{align*}\E d_o(x) &= c\nu W(x,\cdot) = cD_o(x)\\\E d_i(x)&=c\nu W(\cdot,x)=cD_i(x)\end{align*} where $D_o$ and $D_i$ are the normalized mean in- and out-degree functions, c.f. \eqref{eq:Do} and \eqref{eq:Di}.
\end{prop}

 In the following subsections we give properties of the degree function for deterministic and random edge weight functions. We give the results directly in terms of the random measure. These results follow from Theorem~\ref{thm:stc} on STC random measures. 
 
\subsubsection{Deterministic} We give a result on properties of the degree function for deterministic transformations.

\begin{samepage}
\begin{thm}[Deterministic]\label{thm:deter}
Let $N=(\kappa,\nu)$ be a random measure on $(E,\mathscr{E})$ formed by independency $\mathbf{X}=\{X_i:i=1,\dotsb,K\}$. Consider $f\in(\mathscr{E}\otimes\mathscr{E})_{\ge0}$ and define the non-negative function on $(E,\mathscr{E})$ as \begin{equation}\label{eq:degreedeterministic}d(x) = Nf(x,\cdot) = \int_EN(\D y)f(x,y) = \sum_i^K f\circ (x,X_i)\for x\in E\end{equation} Then \begin{enumerate}\item[(i)]
$d(x)$ has Laplace transform \[\varphi_x(\alpha) = \psi(\nu e^{-\alpha f(x,\cdot)})\for\alpha\in\R_{\ge0}\] \item[(ii)] $d(x)$ has mean and variance \begin{align*}\E d(x) &= c\nu f(x,\cdot) \\\Var d(x) &= c\nu(f^2(x,\cdot))+(\delta^2-c)(\nu f(x,\cdot))^2\end{align*} 
\end{enumerate}
\end{thm}
\end{samepage}

\begin{re}[In-degree / out-degree of directed graphs] For asymmetric $f$ resulting in directed graphs, we define $d_o(x)\equiv d(x)=Nf(x,\cdot)$ and $d_i(x)\equiv Nf(\cdot,x)$ as the out-degree and in-degree functions, with Theorem~\ref{thm:deter} holding for $d_i(x)$ by replacing $f(x,\cdot)$ with $f(\cdot,x)$ in expressions. 
\end{re}


\subsubsection{Random} Here we give results on properties of the degree function for random transformations. First we have Bernoulli.

\begin{samepage}
\begin{thm}[Bernoulli]\label{cor:bern} Let $N=(\kappa,\nu)$ be a random measure on $(E,\mathscr{E})$ formed by independency $\mathbf{X}=\{X_i:i=1,\dotsb,K\}$ where $K\sim\kappa$ has pgf $\psi$. Take arbitrary $f\in(\mathscr{E}\otimes\mathscr{E})_{[0,1]}$ and consider the random transformation $\phi:E\times E\mapsto\{0,1\}$ independent of $N$ with marginal transition kernel $Q((x,y),\cdot)=\text{Bernoulli}(f(x,y))$. Take $\phi_x:E\mapsto F$ as $\phi_x(y)=\phi(x,y)$ and define the function \begin{equation}\label{eq:rd}d(x) = (N
\circ\phi_x^{-1})I = \int_EN(\D y)\phi_x(y) = \sum_i^K \phi_x(X_i)\for x\in E\end{equation} Then \begin{enumerate}\item[(i)]$d(x)$ has pgf \[\psi_x(t) =\psi(1-\nu f(x,\cdot)+\nu f(x,\cdot) t )\]  \item[(ii)] $d(x)$ has mean and variance \begin{align*}\E d(x) &= c\nu f(x,\cdot)\\\Var d(x) &= c\nu f(x,\cdot) + (\delta^2-c)(\nu f(x,\cdot))^2\end{align*} 
\end{enumerate}
\end{thm}
\end{samepage}

\begin{re}[In-degree / out-degree of directed graphs]\label{re:indegree} Similar to the deterministic case, for asymmetric $\phi$ of a directed graph, we define $d_o(x)\equiv d(x)=N\phi(x,\cdot)$ and $d_i(x)\equiv N\phi(\cdot,x)$ as the out-degree and in-degree functions. Theorem~\ref{cor:bern} holds for $d_i(x)$ by replacing $f(x,\cdot)$ with $f(\cdot,x)$ in expressions. 
\end{re}

Next, for a digraphon based transformation $\phi$, the degree function is similarly defined, which is unique up to the symmetry of $P$ implied by Lemma~\ref{lem:sym}. 

\begin{samepage}
\begin{thm}[Digraphon]\label{thm:random} Let $N=(\kappa,\nu)$ be a random measure on $(E,\mathscr{E})$ formed by independency $\mathbf{X}=\{X_i:i=1,\dotsb,K\}$ where $K\sim\kappa$ has pgf $\psi$. Let $(F,\mathscr{F})=(\{0,1\}^2,2^F)$ and consider the random transformation $\phi:E\times E\mapsto F$ independent of $N$ with marginal transition kernel $P$ \eqref{eq:P}. Define $\phi_x:E\mapsto F$ as $\phi_x(y)=\phi(x,y)$ and function at $x\in E$ of $f\in\mathscr{F}_{\ge0}$ as \begin{equation}\label{eq:df}d_f(x) = (N\circ\phi^{-1}_x)f = \int_EN(\D y)f(\phi_x(y)) = \sum_i^K f(\phi_x(X_i))\end{equation} Then \begin{enumerate}\item[(i)]$d_f(x)$ has Laplace transform \[\varphi_f^x(\alpha) = \psi(\nu P_xe^{-\alpha f}) = \psi(\int_{E\times F}\nu(\D y)P((x,y),\D z)e^{-\alpha f(z)})\for\alpha\in\R_{\ge0}\] where $P_x=P((x,\cdot),\cdot)$ \item[(ii)] $d_f(x)$ has mean and variance \begin{align*}\E d_f(x) &= c\nu P_x(f)\\\Var d_f(x) &= c\nu P_x(f^2) + (\delta^2-c)(\nu P_x(f))^2\end{align*}
 \end{enumerate} Moreover, these results also hold identically for degree in terms of $\phi^x:E\mapsto F$ defined as $\phi^x(y)=\phi(y,x)$ and functions $f\in\mathscr{F}_{\ge0}$ having either $(0,1),(1,0)\in\supp(f)$ or $(0,1),(1,0)\notin\supp(f)$.\end{thm}
\end{samepage}

We have the following corollary for indicator functions. 
\begin{cor}[pgf] Take $f=\ind{A}$ for $A\subseteq F$. Then $d_f(x)$ has pgf \[\psi_f^x(t) = \psi(1-\nu P_x(A) + \nu P_x(A)t)\] Moreover this representation is unique when either $(0,1),(1,0)\in A$ or $(0,1),(1,0)\notin A$.
\end{cor}\begin{proof}\begin{align*}\psi_f^x(t)&=\psi(\nu P_xt^{\ind{A}})\\&=\psi(\int_{E\times F}\nu(\D y)P((x,y),z)t^{\ind{A}(z)})\\&=\psi(\int_{E\times F}\nu(\D y)\left(P((x,y),z)\ind{A^c}(z) + P((x,y),z)\ind{A}(z)t\right))\\&=\psi(1-\nu P_x(A) + \nu P_x(A)t)\end{align*} 
\end{proof}

Note that for diffuse $\nu$, \[\nu P_x = \sum_{(a,b)\in F}\nu f_{ab}(x,\cdot)\] and for atomic $\nu$,
 \begin{align*}\nu P_x &= \int_{E\times F}\nu(\D y)P((x,y),\D z)\\&=\int_{E\times F}\nu(\D y)(P((x,x),\D z)\ind{}(x=y)+P((x,y),\D z)\ind{}(x\ne y))\\&=\int_{E\times F}\nu(\D y)(R((x,x),\D z)\ind{}(x=y)+Q((x,y),\D z)\ind{}(x\ne y))\\&=\int_{E}\nu(\D y)\ind{}(x=y)\int_FR((x,x),\D z)+\int_{E\times F}\nu(\D y)Q((x,y),\D z)\ind{}(x\ne y)\\&=\nu\{x\} + \sum_{y\in E}\sum_{(a,b)\in F}\nu\{y\}f_{ab}(x,y)\ind{}(x\ne y)\end{align*}

The definition of $f$ encodes the degree function. Take $A=\{(0,1),(1,0)\}$ and $B=\{(1,1)\}$. Then the degree of $f=\ind{A}$ at $x\in E$ is denoted $d_A(x)$ and encodes the number of directed edges of $G(\phi)$ at $x\in E$. Similarly, the degree of $f=\ind{B}$ at $x$, denoted $d_B(x)$, encodes the number of bidirected edges of $G(\phi)$ at $x$. Putting $C=A\cup B$, the degree of $f=\ind{C}$ at $x$, $d_C(x)$, encodes the total number of edges, directed and bidirected, of $G(\phi)$ at $x$.

Lastly, we have degree for the Poisson transformation. Note that Remark~\ref{re:indegree} applies for in- and out-degree functions of directed graphs.

\begin{samepage}
\begin{thm}[Poisson]\label{thm:poisson} Let $N=(\kappa,\nu)$ be a random measure on $(E,\mathscr{E})$ formed by independency $\mathbf{X}=\{X_i:i=1,\dotsb,K\}$ where $K\sim\kappa$ has pgf $\psi$. Take arbitrary $f\in(\mathscr{E}\otimes\mathscr{E})_{\ge0}$ and consider the Poisson random transformation $\phi:E\times E\mapsto\N_{\ge0}$ independent of $N$ with marginal transition kernel $Q((x,y),\cdot)=\text{Poisson}(f(x,y))$. Take $\phi_x:E\mapsto F$ as $\phi_x(y)=\phi(x,y)$ with pgf $\psi_x^f(t)=\exp_-\nu f(x,\cdot)(1-t)$ and define the function \begin{equation}\label{eq:pd}d(x) = (N\circ\phi_x^{-1})I = \int_EN(\D y)\phi_x(y) = \sum_i^K \phi_x(X_i)\for x\in E\end{equation} Then \begin{enumerate}\item[(i)]$d(x)$ has pgf \[\psi_x(t) =\psi\circ\psi_x^f(t)\]  \item[(ii)] $d(x)$ has mean and variance \begin{align*}\E d(x) &= c\nu f(x,\cdot)\\\Var d(x) &= c\nu(f^2(x,\cdot)+f(x,\cdot)) + (\delta^2-c)(\nu f(x,\cdot))^2\end{align*} 
\end{enumerate}
\end{thm}
\end{samepage}

\subsubsection{Degree examples} We give a few examples of degree functions for different graphons of Bernoulli transformations. 

\paragraph{Example: Erd\"{o}s-Renyi} Take $f(x,y)=\ind{}(x\ne y)p$. For Erd\"{o}s-Renyi $G(n,p)=G(\phi)$ with $\kappa=\delta_n$, we have $\E d(x) = np$ and $\Var d(x) = np(1-p)$ for $x\in E$. The degree function $d(x)$ has pgf \[\psi_x(t) = (1-p+pt)^n\] so that $d(x)\sim\text{Binomial}(n,p)$ for $x\in E$. For $x\in\mathbf{X}$, we have $\E d(x) = (n-1)p$, $\Var d(x)=(n-1)p(1-p)$, and $d(x)\sim\text{Binomial}(n-1,p)$. 




\paragraph{Example: Power law} Take $E=[0,1]$, $\nu=\Leb$, and $f(x,y)=(1+bx)^{-2}(1+by)^{-2}$ for $b>0$. Then the degree has mean \[\E d(x) = c\nu f(x,\cdot) = \frac{c}{(1+b)(1+bx)^2}\sim x^{-2}\] and  variance \[\Var d(x) = c\nu f^2(x,\cdot) + (\delta^2-c)(\nu f(x,\cdot))^2 = \frac{c b^2+3 (1+b) \delta^2}{3 (1+b)^3 (1+b x)^4}\sim x^{-4}\] The mean total edge weight is given by \begin{align*}\E Nd &=c(\nu\times I)f + (c^2+\delta^2-c)(\nu\times\nu)f \\&= c\frac{3+b(3+b)}{3(1+b)^3} + \frac{c^2+\delta^2-c}{(1+b)^2}\\&=\frac{cb^2+3c^2(1+b)+3(1+b)\delta^2}{3(1+b)^3}\end{align*} 

\paragraph{Example: Power law (continuation)} Define $Q((x,y),\cdot)=\text{Bernoulli}(f(x,y))$. Then $d$ has the same mean \[\E d(x) = c\nu f(x,\cdot) = \frac{c}{(1+b)(1+bx)^2}\sim x^{-2}\]  and the variance is given by\[\Var d(x) = c\nu f(x,\cdot) + (\delta^2-c)(\nu f(x,\cdot))^2=\frac{c b (1+(1+b)(2+b x)x)+\delta^2}{(1+b)^2 (1+b x)^4}\sim x^{-2}\] The mean number of edges is given by \[\E (M\circ\phi^{-1})I = c(\nu\times I)f + (c^2+\delta^2-c)(\nu\times\nu)f = \frac{cb^2+3c^2(1+b)+3(1+b)\delta^2}{3(1+b)^3}\] which is also the same. 

\subsection{Degree distribution}\label{sec:degreedis} An important quantity is the distribution of the degree function for random inputs, the degree distribution (DD). We give the pgf of the degree distribution (DD pgf) for various transformations. 

 The following theorem gives the DD pgf for deterministic and Bernoulli random transformations.
 
\begin{thm}[Bernoulli DD]\label{prop:dis}Consider deterministic $f\in(\mathscr{E}\otimes\mathscr{E})_{[0,1]}$ or random $\phi:E\times E\mapsto\{0,1\}$ with $\phi(x,y)\sim\text{Bernoulli}(f(x,y))$ for $f\in(\mathscr{E}\otimes\mathscr{E})_{[0,1]}$. Let $Y=d(X)$ for degree function $d$ \ref{eq:rd} and $X\sim\nu$. Then the pgf of $Y$ is given by \[\psi_Y(t)=\int_E\nu(\D x)\psi(1-\nu f(x,\cdot)+\nu f(x,\cdot) t)\] \end{thm}

\begin{re}[Directed] If $f$ is asymmetric, then we denote $Y_o=d(X)$ as the out-degree of $X$, with pgf $\psi_{Y_o}(t)=\int_E\nu(\D x)\psi(1-\nu f(x,\cdot)+\nu f(x,\cdot) t)$, and $Y_i=d_i(X)$ as the in-degree of $X$, with pgf  $\psi_{Y_i}(t)=\int_E\nu(\D x)\psi(1-\nu f(\cdot,x)+\nu f(\cdot,x) t)$
\end{re}

Next we have the DD pgf for the diagraphon transformation. 

\begin{thm}[Digraphon DD] For digraphon random transformation $\phi:E\times E\mapsto F$ where $(F,\mathscr{F})=(\{0,1\}^2,2^F)$ and indicator function $f=\ind{A}$, $A\subseteq F$, the pgf of $Y=d_f(X)$ for degree function $d_f$ \eqref{eq:df} is given by \[\psi_Y(t)=\int_E\nu(\D x)\psi(1-\nu P_x(A)+\nu P_x(A)t)\] This is unique when $(0,1),(1,0)\in A$ or $(0,1),(1,0)\notin A$.
\end{thm}

Lastly we give the DD pgf for the Poisson transformation. 

\begin{thm}[Poisson DD]\label{prop:poisson}Consider the Poisson transformation $\phi:E\times E\mapsto\N_{\ge0}$ with $\phi(x,y)\sim\text{Poisson}(f(x,y))$ for $f\in(\mathscr{E}\otimes\mathscr{E})_{\ge0}$. Let $Y=d(X)$ for degree function $d$ \eqref{eq:pd} and $X\sim\nu$. Then the pgf of $Y$ is given by \[\psi_Y(t)=\int_E\nu(\D x)\psi(\exp_-\nu f(x,\cdot)(1-t))\] \end{thm}

Recall the probabilities are extracted using the Cauchy formula. 

\begin{re}[Coefficients] By Theorem~\ref{thm:dis}, the probability $\P(Y=k)$ may be attained through the Cauchy formula as  \[\P(Y=k)=\frac{1}{2\pi i}\oint_C\frac{\psi_Y(t)}{t^{k+1}}\D t\for k\in\N_{\ge0}\]  
\end{re}

\paragraph{Example: Bernoulli transformation, Poisson random measure}Consider $\kappa=\text{Poisson}(c)$ and Bernoulli $\phi$. Then $Y=d(X)$ for $X\sim\nu$ has pgf \[\psi_Y(t) =\int_E\nu(\D x)\exp_-c\nu f(x,\cdot)(1-t)\] and mean and variance \begin{align*}\E Y &= \int_E\nu(\D x)c\nu f(x,\cdot) = c(\nu\times\nu)f\\\Var Y &= \int_E\nu(\D x) c\nu f(x,\cdot)(1+c\nu f(x,\cdot)) - (\E Y)^2\\&=\E Y + c^2\int_E\nu(\D x)(\nu f(x,\cdot))^2 - (\E Y)^2\end{align*} In some cases $\psi_Y$ may be closed-form. For example, let $f(x,y)=g(x)g(y)$, so that  \[\psi_Y(t) =\int_E\nu(\D x)\exp_-c\nu(g)g(x) (1-t)\] Take $E=[0,1]$, $\nu=\Leb$, and power-law $g(x)=(1+bx)^{-2}$, $b>0$. With $\nu(g)=(1+b)^{-1}$, the pgf is \[\psi_Y(t) = \frac{1}{b}\left(\sqrt{\pi } \sqrt{c \nu(g) (1-t)} \left(\text{erf}\left(\frac{\sqrt{c \nu(g) (1-t)}}{1+b}\right)-\text{erf}\left(\sqrt{c \nu(g) (1-t)}\right)\right)+(1+b) e^{-\frac{c \nu(g) (1-t)}{(1+b)^2}}-e^{-c \nu(g) (1-t)}\right)\]

\subsection{Giant components}\label{sec:giant} A giant component of a graph is a connected component that contains a finite fraction of the graph's vertices. 

\begin{thm}[GC] Consider an undirected random graph $G(\phi)$ with degree function $d$. Put $Y=d(X)$ for $X\sim\nu$. Then $G(f)$ contains a giant component almost surely iff $\E Y^2-2\E Y>0$. \end{thm}\begin{proof}\citep{gc}. \end{proof}


\paragraph{Example: Poisson random measure} For example,  consider $\kappa=\text{Poisson}(c)$. Then the condition for a GC is \[c\int_E\nu(\D x) (\nu f(x,\cdot))^2  - (\nu\times\nu)f>0\] For $f(x,y)=p$, this reduces to \[c>1,\quad \frac{1}{c}<p\le 1\] For $\kappa=\text{Dirac}(n)$, the relation is \[(n-1)\int_E\nu(\D x) (\nu f(x,\cdot))^2-(\nu\times\nu)f>0\] and for $f(x,y)=p$, reduces to \[n\ge3,\quad \frac{1}{n-1}<p\le 1\]

\subsection{Active vertices}\label{sec:part} An important quantity for the random graph $G(f)$ are the vertices appearing in its support, \emph{active vertices}, that is, the vertices of those (active) edges where $f>0$. We give the mean number of active vertices in $G(f)$ for deterministic and random $f$. The results are similar. Note these these results hold also for directed graphs (following from asymmetric $f$) as out-degree active vertices; for in-degree active vertices, the results hold when replacing $f(x,\cdot)$ by $f(\cdot,x)$ in expressions. 

\begin{prop}[Active vertices]\label{prop:participate} Let $(M=N\times N,\ind{A},I)$ form the random graph $G=G(\ind{A})$ for arbitrary $A\subseteq E\times E$. Form $d(x)=N\ind{A}(x,\cdot)=K_A^x$, define $g(x)=\ind{\ge1}(K_A^x)$, and put $a(x)=\nu(\ind{A}(x,\cdot))$. Then the number of active vertices $v(G)=Ng$ of $G$ has mean \[\E v(G)=c\int_E\nu(\D x)(1-\psi(1-a(x))(1-\ind{A}(x,x)))\]
\end{prop}
\begin{proof}
\begin{align*}\E Ng &=c\int_{E}\nu(\D x)\P(K_A^x\ge1|f(x,x))\\&=c\int_E\nu(\D x)(1-\P(K_A^x=0|f(x,x))\\&=c\int_E\nu(\D x)(1-\psi_A^x(0)(1-\ind{A}(x,x)))\end{align*} 
\end{proof}

For Poisson $\kappa$, making use of $a(x)=\nu(\ind{A}(x,\cdot))$, this is \[\E v(G) = c\nu(1-e^{-ca}) + c\nu e^{-ca}\ind{A}(\cdot,\cdot)\] To avoid self-edges, we can define $A$ through considering subspaces $B\subseteq E\times E$ and $C=\{(x,y)\in E\times E: x\ne y\}$ and setting $A= B\cap C$, so that $\E M(A) = (c^2+\delta^2-c)(\nu\times\nu)(A)$ and $\E v(G) = c\nu(1-\psi_A^\bullet(0))$. In this case, Poisson has mean $\E v(G) = c\nu(1-e^{-ca})$. 
 
Now consider random indicators through the Bernoulli transformation. The following result gives the mean number of vertices active in the graph edges. 

\begin{prop}[Active vertices]\label{prop:participate} Let $(M=N\times N,\phi,I)$ form the random graph $G=G(\phi)$ with  Bernoulli transformation $\phi$ based on $f\in(\mathscr{E}\otimes\mathscr{E})_{[0,1]}$. Form $d(x)=N\phi(x,\cdot)=K_x$ and define $g(x)=\ind{\ge1}(K_x)$. Then the number of active vertices $v(G)=Ng$ of $G$ has mean \[\E v(G) =c\int_E\nu(\D x)(1-\psi(1-\nu f(x,\cdot))(1-f(x,x)))\]
\end{prop}

Consider random graph $G(\phi)$ with Bernoulli transformation $\phi$ based on graphon $f$. For Poisson, the mean of $v(G)$, putting $a(x)=\nu f(x,\cdot)$, is given by \[\E v(G) = c\nu(1-e^{-ca}) + c\nu e^{-ca}f(\cdot,\cdot)\] 

\paragraph{Example: Power-law graphon} Let $E=[0,1]$ and $\nu=\Leb$.  Consider Poisson and graphon $f(x,y)=\ind{}(x\ne y)g(x)g(y)$ for $g(x)=(1+bx)^{-2}$ where $\nu(g)=(1+b)^{-1}$. Put $a(x)=\nu(g)g(x)$. We have  \[\nu e^{-ca} = \frac{1}{b}\left(\sqrt{\pi } \sqrt{c \nu(g)} \left(\text{erf}\left(\frac{\sqrt{c \nu(g)}}{1+b}\right)-\text{erf}\left(\sqrt{c \nu(g)}\right)\right)+(1+b) e^{-\frac{c \nu(g)}{(1+b)^2}}-e^{-c \nu(g)}\right)\xrightarrow{b\rightarrow\infty}1\] so that $\E v(G) = c\nu(1-e^{-ca})\xrightarrow{b\rightarrow\infty}0$. Note that \[\nu e^{-ca}\sim \exp_-c \nu(g)/(1+b)^2 = \exp_-c/(1+b)^3\quad\text{as}\quad b\rightarrow\infty\]

\paragraph{Example: Exponential graphon} Let $E=[0,1]$ and $\nu=\Leb$.  Consider Poisson and graphon $f(x,y)=\ind{}(x\ne y)g(x)g(y)$ for $g(x)=\exp_-bx$, $b>0$, where $\nu(g)=(1-e^{-b})/b$ (which coincidentally is the Laplace transform of $\nu$ in $b$). Put $a(x)=\nu(g)g(x)$. Then we have \[\nu e^{-ca} = \frac{1}{b}\left(\text{Ei}(-c \nu(g))-\text{Ei}\left(-c e^{-b} \nu(g)\right)\right)\xrightarrow{b\rightarrow\infty}1\] where $\text{Ei}$ is the exponential integral. Hence $\E v(G) = c\nu(1-e^{-ca})\xrightarrow{b\rightarrow\infty}0$. We have \[\nu e^{-ca}\sim\exp_-c\nu(g)e^{-b}=\exp_-ce^{-b}(1-e^{-b})/b\sim\exp_-ce^{-b}/b\quad\text{as}\quad b\rightarrow\infty\] Because $1/(1+b)^3$ dominates $e^{-b}/b$ for large $b$, the exponential graphon  (asymptotically) has fewer active vertices than power-law. 

\paragraph{Example: Erd\"{o}s-Renyi graphon} Let $E=[0,1]$ and $\nu=\Leb$.  Consider Poisson and graphon $f(x,y)=\ind{}(x\ne y)p$. Put $a(x)=p$ so that \[\nu e^{-ca} = e^{-cp}\xrightarrow{p\rightarrow0}1\] Hence $\E v(G) = c\nu(1-e^{-ca})=c(1-e^{-cp})\xrightarrow{p\rightarrow0}0$.

\subsection{Empirical graphs}\label{sec:emp} Consider an observed collection of variates $\mathbf{X}=\{X_1,\dotsb,X_n\}$ taking values in $(E,\mathscr{E})$. We refer to $G(g\circ f)=(\mathbf{X},\mathbf{X}\times\mathbf{X},g\circ f)$ for $g\circ f\in(\mathscr{E}\otimes\mathscr{E})_{\ge0}$ as an \emph{empirical} graph. 

Define the empirical distribution (probability measure) $F_n$ on $(E,\mathscr{E})$ corresponding to the $X_1,\dotsb,X_n$ as \begin{equation}\label{eq:emp}F_n(A)= \frac{1}{n}\sum_i^n\ind{A}(X_i)\for A\in\mathscr{E}\end{equation} The collection $\mathbf{X}$ forms the deterministic counting measure $N$ on $(E,\mathscr{E})$ as $N=nF_n$, and the collection $\mathbf{X}\times\mathbf{X}$ forms the deterministic counting measure $M=N\times N$ on $(E\times E,\mathscr{E}\otimes\mathscr{E})$ as $M=n^2(F_n\times F_n)$. For random $f$ with marginal transition kernel $Q$, we have \[\E(M\circ f^{-1})g = \int_{E\times E}N(\D x)N(\D y)\E (g\circ f)(x,y) = n^2(F_n\times F_n)Qg\for g\in\mathscr{F}_{\ge0}\]

 \paragraph{Example: Soft fixed-degree graphs} The empirical graph can be used to represent random graphs having fixed mean degree sequence. Consider observed fixed degree sequence $\mathbf{D}=\{D_1,\dotsb,D_n\}$ with total degree $m=D_1+\dotsb+D_n$. Let $N$ be the deterministic counting measure $N$ formed by $\mathbf{D}$, and let $M=N\times N$ be the deterministic product counting measure formed by $\mathbf{D}\times\mathbf{D}$. Let $G(\phi)=(\mathbf{D},\mathbf{D}\times\mathbf{D},\phi)$ be an empirical random graph with Poisson transformation $\phi(x,y)\sim\text{Poisson}(f(x,y))$ for $f\in(\mathscr{E}\otimes\mathscr{E})_{\ge0}$ as $f(x,y)=xy/m$. This $f$ preserves the total $\E M\phi=m$ and individual degrees $\E N\phi(x,\cdot)=x$ in the mean. Thus $G(\phi)$ generates random graphs with mean fixed degree sequence $\mathbf{D}$. We refer to Poisson $\phi$ with such $f$ as a \emph{soft} fixed-degree transformation, where soft indicates the graphs respect the degree sequence in the mean. 
 
Now we consider $\mathbf{X}=\{X_1,\dotsb,X_n\}$ to be an independency of identically $\nu$-distributed $E$-valued random variables forming empirical distribution $F_n$. Let $N=(\kappa,F_n)$ on $(E,\mathscr{E})$ be STC random measure formed by $\mathbf{X}'=\{X_i':i=1,\dotsb,K\}$ of $K\sim\kappa=\text{Dirac}(c)$ iid $X'_i\sim F_n$. We have \[\E Ng = cF_n(g)\xrightarrow{n\rightarrow\infty}c\nu g\quad\text{almost surely}\for g\in\mathscr{E}_{\ge0}\] 

The degree function $d(x)=Nf(x,\cdot)$ for deterministic $f$ has mean and variance \begin{align*}\E d(x) &= c F_nf(x,\cdot)\\\Var d(x) &= c\Var_nf(x,\cdot)\end{align*} For degree function $d(x)=N\phi(x,\cdot)$ for Bernoulli $\phi$ with graphon $f\in(\mathscr{E}\otimes\mathscr{E})_{[0,1]}$, we have \begin{align*}\E d(x) &= c F_nf(x,\cdot)\\\Var d(x) &= cF_nf(x,\cdot)(1-F_nf(x,\cdot))\end{align*} Here $d(x)=d_o(x)$ is the out-degree function. The results hold for the in-degree function $d_i(x)$ by replacing $f(x,\cdot)$ with $f(\cdot,x)$ in expressions. 

 Let $M=N\times N$ be the product random measure formed by $\mathbf{X}'\times\mathbf{X}'$. We refer to $G(f)=(\mathbf{X}',\mathbf{X}'\times\mathbf{X}',f)$ as a \emph{sampled empirical} graph. The mean of $Mf$ and its convergence are given by \[\E Mf = c\,(F_n\times I)f + c(c-1)\,(F_n\times F_n)f\,\,\,\xrightarrow{n\rightarrow\infty}\,\,\,c\,(\nu\times I)f + c(c-1)\,(\nu\times\nu)f\quad\text{almost surely}\] 
 
 Oftentimes $c=n$, in which case we refer to a \emph{bootstrap} sampled empirical graph, such that $Mf$ has mean \[\E Mf = n\,(F_n\times I)f + n(n-1)\,(F_n\times F_n)f = \sum_i^nf(X_i,X_i) + \frac{n-1}{n}\sum_{i}^n\sum_j^nf(X_i,X_j)\] 

 \subsection{Directed acyclic graphs}\label{sec:dag}
 
 A \emph{path} in a directed graph is a sequence of edges where the end vertex of each edge in the sequence is equal to the starting vertex of the next edge in the sequence. Such a path is a \emph{directed cycle} if the starting vertex of the first path edge is equal to the ending vertex of the last path edge. A \emph{directed acyclic graph} (DAG) is a directed graph having no directed cycles. A graph is a DAG iff its adjacency matrix is lower-triangular with zero diagonal.
 
Consider the STC random measure $N=(\kappa,\nu)$ on $(E,\mathscr{E})$ formed by independency $\mathbf{X}$ and its product $M=N\times N$ on $(E\times E,\mathscr{E}\otimes\mathscr{E})$ formed by $\mathbf{X}\times\mathbf{X}$. Assume $E$ is ordered. Consider subspace $A=\{(x,y)\in E\times E: x<y\}\subset E\times E$. Let $(M,\phi,I)$ form the STC random graph $G=G(\phi)=(\mathbf{X},\mathbf{X}\times\mathbf{X},\phi)$ where $\phi$ is a Bernoulli transformation with restricted graphon $f_A=f\ind{A}$ of graphon $f$, forming transition kernel $Q$, i.e. $\phi(x,y)\sim Q((x,y),\cdot)=\text{Bernoulli}(f(x,y)\ind{A}(x,y))$. $G$ is a DAG. The number of edges of $G$ is $e(G)=(M\circ\phi^{-1})I$ with mean \[\E e(G) = (c^2+\delta^2-c)(\nu\times\nu)f_A\]  The out-degree and in-degree functions in $x\in E$ are given by $d_o(x) = N\phi(x,\cdot)$ and $d_i(x)=N\phi(\cdot,x)$ respectively, with means $\E d_o(x) = c\nu f_A(x,\cdot)$ and $\E d_i(x) = c\nu f_A(\cdot,x)$. Hence the underlying graphon $f$ encodes the properties of $G$. In many applications, sparse DAGs are demanded, which are satisfied using suitable sparsifying graphons. For separable graphon $f=g\times g$ and putting $a_x(y)=\ind{A}(x,y)$ and $b_x(y)=a_y(x)$, the mean out and in degrees are $\E d_o(x)=c\nu(ga_x)g(x)$ and $\E d_i(x)=c\nu(gb_x)g(x)$ respectively.  

\subsection{Rewired graphs}\label{sec:resample} We describe an iterative resampling procedure (\emph{resampler}) that rewires Bernoulli random graphs and preserves the product mean measure across iterations for fixed graphon. We define a rewiring transform that completely randomly rewires a fixed number of random vertices at each iteration. 

First we define the initial condition for the resampler. It is the random graph $G=G(\phi)=(\mathbf{X},\mathbf{X}\times\mathbf{X},\phi)$ formed by $(M=N\times N,\phi,I)$ with $N=(\kappa,\nu)$ on $(E,\mathscr{E})$ and random $\phi(x,y)\sim\text{Bernoulli}(f(x,y))$ for graphon $f$. The adjacency matrix is $\mathbf{A}=\phi(\mathbf{X}\times\mathbf{X})\in\{0,1\}^{K\times K}$. 

Each successive iteration of the resampler is defined as follows. Let $J\subseteq\{1,\dotsb,K\}$ be a random subset of size $n\in\{1,\dotsb,K\}$. Let $\mathbf{X}=\{X_i:i=1,\dotsb,K\}$ be the previous iteration's independency, and define resampled variates for the current iteration as \begin{equation}X^*_i = \begin{cases}Z_i\sim\nu&\text{if }i\in J\\X_i&\text{otherwise}\end{cases}\for i\in\{1,\dotsb,K\}\end{equation} These form the independency of $n$-resampled variates $\mathbf{X}^*=\{X_i^*:i=1,\dotsb,K\}$. Essentially this is replacement of certain variates with fresh random realizations while freezing the remaining. Given the previous iteration's transform $\phi$ and the current iteration's resampled independency $\mathbf{X}^*\times\mathbf{X}^*$, define the \emph{resampling transform} $\phi^*$ as \begin{equation}\phi^*(X_i^*,X_j^*) = \begin{cases}Z_{ij}\sim\text{Bernoulli}(f(X_i^*,X_j^*))&\text{if }i\in J\text{ and/or }j\in J\\\phi(X_i^*,X_j^*)&\text{otherwise}\end{cases}\for i,j=1,\dotsb,K \end{equation} such that $\mathbf{A}^*=\phi^*(\mathbf{X}^*\times\mathbf{X}^*)$ is the $n$-resampled adjacency matrix. Note that if $\phi$ is symmetric, then $\phi^*$ is taken to be symmetric. Finally we define the \emph{($n$-)rewiring transform} $S_n$ of STC random graph $G(\phi)=(\mathbf{X},\mathbf{X}\times\mathbf{X},\phi)$ as the mapping \begin{equation}\label{eq:sn}S_n\circ (\mathbf{X},\mathbf{X}\times\mathbf{X},\phi) \equiv (\mathbf{X}^*,\mathbf{X}^*\times\mathbf{X}^*,\phi^*)\end{equation} and denote $G^*=G(\phi^*)=S_n\circ G(\phi)$ as the \emph{$n$-rewired graph} of $G(\phi)$. By construction, the mean adjacency matrix is unchanged: $\E\mathbf{A}^*=\E\mathbf{A}$. Hence, the sequence of rewired graphs $(S_n^i(G(\phi)))_{i\ge0}$ formed by iterating, where $S_n^i$ is the $i$-fold composition of $S_n$ with itself, preserves the mean measure. 

 Note that other rewiring transforms can be defined. For example, we can define \begin{equation}\phi^*(X_i^*,X_j^*) = \begin{cases}Z_{ij}\sim\text{Bernoulli}(f(X_i^*,X_j^*))&\text{if }i,j\in J\\\phi(X_i^*,X_j^*)&\text{otherwise}\end{cases}\for i,j=1,\dotsb,K \end{equation} which resamples only the edges in $J\times J$ at each iteration.

\section{FAIW random measure graph models}\label{sec:farwg}

Another class of random graphs are those defined in terms of products of random counting measures having fixed atoms and random integer weights (FAIW). This class is equivalent to Gibbs random fields on lattices (Remark~\ref{re:part}). It is prototypical to unlabeled graphs. 
\begin{samepage}
\begin{de}[FAIW random measure graph]\label{def:rg} Consider the FAIW random counting measure $N$ on $(E,\mathscr{E})$ with countable set of fixed atoms $D\subseteq E$ and independent non-negative integer weights $W=\{W_x: x\in D\}$ distributed $W_x\sim\kappa_x$ with mean $c_x$ and variance $\delta_x^2$. Now consider the product random counting measure $M=N\times N$ on $(E\times E,\mathscr{E}\otimes\mathscr{E})$. For functions $g\circ f$ and $h_W$ in $(\mathscr{E}\otimes\mathscr{E})_{\ge0}$ where $h_W(x,y)=W_xW_yg\circ f(x,y)$, the triple $(M,f,g)$ defines a non-negative weighted \emph{random graph} $G_W=G(h_W)=(D,D\times D,h_W)$ where $D$ is the vertex set, $D\times D$ is the edge set, and $h_W$ is the edge weight function. 
\end{de}
\end{samepage}

 Given FAIW graph $G=G(h_W)$, similar definitions apply as with STC graphs. Here we give only the basic properties of FAIW graphs (the remaining constructs may be ported). 
 
 The support of $G(h_W)$ is $\supp(G(h_W)) = \{h_W>0\} = \{(x,y)\in D\times D: h_W(x,y)>0\}$. The number of active edges of $G$ is given by $e(G) \equiv |\supp(G)| = (M\circ f^{-1})\ind{\supp(g)}$, and the weight of $G$ is $w(G) \equiv (M\circ f^{-1})g$. If $g\circ f$ is symmetric, then $G(g\circ f)$ is symmetric, and edge count and weight are defined relative to the normalized product FAIW random measure $M^*$ of $M$ (Remark~\ref{re:norm2}) as $e^*(G)\equiv(M^*\circ f^{-1})\ind{\supp(g)}$ and $w^*(G)\equiv(M^*\circ f^{-1})g$ respectively.
 
 Consider subset $A\subseteq D$ and the restricted measures $N_A(\cdot)=N(A\cap\cdot)$ and $M_{A\times A}=N_A\times N_A$. Then the triple $(M_{A\times A},f,g)$ forms the random graph $G_{A\times A}=G_{A\times A}(h)=(D\cap A,(D\times D)\cap(A\times A),h)$ as the restriction $G_{A\times A}$ of the graph $G$ to subset $A\times A\subseteq D\times D$. 
 
The remainder of this section is organized as follows. We discuss deterministic and random transformations (\ref{sec:tran}), degree functions (\ref{sec:fiawdegree}), fixed-degree graphs (\ref{sec:fixed}), and Bernoulli thinning of adjacency matrices (\ref{sec:thin}).

\subsection{Transformations}\label{sec:tran} As with STC random graphs, we discuss deterministic and random transformations. 

\subsubsection{Deterministic} When $f$ is deterministic, the mean edge count and weight of $G$ are \[\E(M\circ f^{-1})\ind{\supp(g)} = \sum_{(x,y)\in D\times D}Z_{xy}\ind{\supp(g)}\circ f(x,y) = \sum_{(x,y)\in D\times D}Z_{xy}\ind{}(g\circ f(x,y) > 0)\] and \[\E(M\circ f^{-1})g = \sum_{(x,y)\in D\times D}Z_{xy}\, g\circ f(x,y)\] respectively where $Z_{xy} = \E W_xW_y$.

\subsubsection{Random} Letting random $\phi$ have marginal transition kernel $Q$, and defining measure $Q_{xy}(\cdot)=Q((x,y),\cdot)$, the mean edge count and weights are respectively given by \[ \E(M\circ \phi^{-1})\ind{\supp(g)} = \sum_{(x,y)\in D\times D}Z_{xy}Q_{xy}(\supp(g))\] and \[\E(M\circ \phi^{-1})g = \sum_{(x,y)\in D\times D}Z_{xy}Q_{xy}(g)\] 


We discuss some random transformations for product FAIW random measures. 

\paragraph{Bernoulli} Consider the graphon $f\in(\mathscr{E}\otimes\mathscr{E})_{[0,1]}$ based Bernoulli transformation $\phi:D\times D\mapsto\{0,1\}$ with marginal Bernoulli kernel $Q((x,y),\cdot)\sim\text{Bernoulli}(f(x,y))$. Then  \[ \E(M\circ \phi^{-1})I =  \sum_{(x,y)\in D\times D}Z_{xy}f(x,y)\]
 
\paragraph{Poisson} Another case is the Poisson transformation $\phi:D\times D\mapsto\N_{\ge0}$ with marginal Poisson kernel $Q((x,y),\cdot)\sim\text{Poisson}(f(x,y))$ for $f\in(\mathscr{E}\otimes\mathscr{E})_{\ge0}$. Then  \[ \E(M\circ \phi^{-1})I =  \sum_{(x,y)\in D\times D}Z_{xy}f(x,y)\]

\paragraph{Exponential} Let $\phi:D\times D\mapsto \{0,1\}$ be a random transformation specified through the exponential family of random graphs \citep{exp}. In particular, the joint law of $\phi=\{\phi xy: (x,y)\in D\times D\}$ is given as \[\P(\phi=A|\theta) = \frac{e^{\theta\cdot s(A)}}{Z(\theta)} \for A\in\{0,1\}^{|D|\times |D|}\] where $s(A)$ is a vector of sufficient statistics of the adjacency matrix $A$, $\theta$ is a vector of parameters corresponding to the sufficient statistics, and $Z(\theta)=\sum_{A}e^{\theta\cdot s(A)}$ is the normalizing constant, a partition function. In practice $Z(\theta)$ is infeasible to compute and random realizations of $\phi$ are generated using the Metropolis-Hastings algorithm. The marginal probabilities are also intractable and thus the marginal transition kernel $Q$ of $\phi$ cannot be meaningfully defined.

\paragraph{Preferential attachment}  Let $\phi:D\times D\mapsto \{0,1\}$ be a random transformation specified through a preferential attachment mechanism, e.g.,  Barb\'{a}si-Albert \citep{pref}, where the law of $\phi=\{\phi xy: (x,y)\in D\times D\}$ is implied through the iterative generator. This iteration process is initialized with the first, say $n$, vertices of $D$ and some edge set (perhaps empty). For each successive element of $D$, its space of edges is randomly sampled, materializing the corresponding row and column of the adjacency array. This process is repeated until $D$ is exhausted. 

\paragraph{Multinomial} Let $\phi_n:D\times D\mapsto\N_{\ge0}$ be a multinomial transformation: $\phi_n=\{\phi_n xy: (x,y)\in D\times D\}$ has law \[\P(\phi_n xy = a_{xy}: (x,y)\in D\times D) = n!\prod_{(x,y)\in D\times D}\frac{p_{xy}^{a_{xy}}}{a_{xy}!}\] where $(p_{xy})$ are probabilities $\sum_{x,y}p_{xy}=1$ and $a_{xy}\in\{0,1,\dotsb,n\}$. $\phi_n$ has a binomial marginal transition kernel: $\phi_nxy\sim Q_n((x,y),\cdot)=\text{Binomial}(n,p_{xy})$. Then  \[ \E(M\circ \phi_n^{-1})I =  n\sum_{(x,y)\in D\times D}Z_{xy}p_{xy}\] This transformation is asymmetric and yields directed multigraphs. A symmetric transformation yielding undirected multigraphs without loops is defined as follows: $n$ is even such that a multinomial transformation $\phi_{n/2}$ is defined relative to a triangular of $D\times D$ (not including diagonal) and the result symmetrized to yield $\phi_n$ with $p_{xy}=p_{yx}$. If loops are allowed, then the same procedure is applied to a triangular including the diagonal and the result symmetrized, yielding an undirected multigraph with loops having number of edges bounded by $n$. Essentially the inclusion of loops for undirected graphs destroys the control over the number of edges by the multinomial transform yet does not alter the mean number of edges.



 \subsection{Degree}\label{sec:fiawdegree} Consider FAIW graph $G=G(h_W)=(D,D\times D,h_W)$ formed from $(M,f,g)$. The degree function of $G$ is defined as \[d(x)=N(g\circ f(x,\cdot)) = \sum_{z\in D}W_z\,g\circ f(x,z)\for x\in E\]  For subset $A\subseteq D$ and restricted random measures $N_A(\cdot)=N(A\cap\cdot)$ and $M_{A\times A}=N_A\times N_A$, the degree function of $G$ restricted to $A\times A$ is \[d_A(x) = \sum_{z\in D\cap A}\ind{A}(x)\,W_z\,g\circ f(x,z)\for x\in E\]
 
\subsubsection{Deterministic} Consider deterministic $f:E\times E\mapsto\R_{\ge0}$ and $g=I$.

 The degree function has mean \[\E d(x) = \sum_{z\in D}c_zf(x,z)\for x\in E\] and Laplace transform \[\varphi_{x}(\alpha) =\prod_{z\in D}\varphi_{W_z}(\alpha f(x,z)) \for\alpha\in\R_{\ge0},\quad x\in E\] 

The degree function $d_A$ of the graph restriction to subset $A\times A\subseteq D\times D$ has mean \[\E d_A(x) = \sum_{z\in D\cap A}c_z\ind{A}(x)f(x,z)\for x\in E\] and Laplace transform \[\varphi_{A}^x(\alpha) =\prod_{z\in D\cap A}\varphi_{W_z}(\alpha \ind{A}(x)f(x,z)) \for\alpha\in\R_{\ge0},\quad x\in E\] 

 
\subsubsection{Random} Consider random $f$ having marginal transition kernel with measure $Q_{xy}(\cdot)=Q((x,y),\cdot)$ on $E\times E$ with support $D\times D$.

 The degree function has mean \[\E d(x) = \sum_{z\in D}c_zQ_{xz}(g)\for x\in E\] Defining the Laplace transform of $f(x,y)$ as \[\varphi_{f(x,y)}(\alpha) = \E e^{-\alpha f(x,y)} = Q_{xy}e^{-\alpha I} = \int_F Q((x,y),\D z)e^{-\alpha z}\for \alpha\in\R_{\ge0}\] the degree function has Laplace transform \[\varphi_x(\alpha) = \prod_{z\in D}(\varphi_{W_z}\circledast\varphi_{f(x,z)})(\alpha) \for\alpha\in\R_{\ge0},\quad x\in E\] where $\circledast$ is the convolution operator.
 
 The restricted degree $d_A$ has mean \[\E d_A(x) = \sum_{z\in D\cap A}c_z\ind{A}(x)Q_{xz}(g)\for x\in E\] and  Laplace transform \[\varphi_A^x(\alpha) = \prod_{z\in D\cap A}(\varphi_{W_z}\circledast\varphi_{\ind{A}(x)f(x,z)})(\alpha) \for\alpha\in\R_{\ge0},\quad x\in E\]

\paragraph{Example: Bernoulli} For Bernoulli $f$ with graphon $W$ and $\kappa_x=\text{Bernoulli}(p_x)$, we have degree mean and pgf  \[\E d(x) = \sum_{z\in D}p_zW(x,z)\for x\in E\]  and \[\psi_x(t) = \prod_{z\in D}(1-p_zW(x,z)+p_zW(x,z)t)\for x\in E\] Hence degree is a limiting (infinite dimensional) Poisson-binomial distribution. 

\subsection{Fixed edge count graphs}\label{sec:fixedge} 

Consider random multigraphs $G$ having fixed edge count $n$. We can represent such random graphs using a FAIW counting measure $N$ on $(E,\mathscr{E})=(\N_{\ge1},2^{\N_{\ge1}})$ with fixed atoms $E$ and unit weights ($\kappa_x=\delta_1$). Let $M=N\times N$ be the FAIW product measure on $(E\times E,\mathscr{E}\otimes\mathscr{E})$ and let $\phi_n:E\times E\mapsto\N_{\ge0}$ be a multinomial transform ($n$ even if symmetric). Then $G(\phi_n)=(E,E\times E,\phi_n)$ is a random multigraph with edge count $n$.
 
\subsection{Fixed degree graphs}\label{sec:fixed} 
Here we discuss undirected random graphs $G$ having fixed \emph{degree sequence} $\mathbf{D}=\{D_1,\dotsb,D_n\}$ where the $D_i\in\N_{\ge0}$. To represent such graphs, we use a FAIW counting measure $N$ on $(E,\mathscr{E})=(\{1,\dotsb,n\},2^{\{1,\dotsb,n\}})$ with fixed atoms $E$ and unit weights $W^1=\{1:x\in E\}$, i.e. $\kappa_x=\delta_1$. Let $M=N\times N$ on $(E\times E,\mathscr{E}\otimes\mathscr{E})$ be the FAIW product measure. 

In particular we form the random graph $G(\phi)=(E,E\times E,\phi)$ from the triple $(M,\phi,I)$, where $\phi:E\times E\mapsto F\subseteq\N_{\ge0}$ is a symmetric random transformation defined such that $d(X_i)=N\phi(X_i,\cdot)=D_i$ for $X_i\in\mathbf{X}$, i.e. a \emph{hard} fixed-degree transformation, where hard means the resultant graph has fixed-degree $\mathbf{D}$ with probability one. The determination of existence of a hard transformation $\phi$ for some $F$ and $\mathbf{D}$ is the same as the matter of the existence of an adjacency matrix $\mathbf{A}=(\phi(i,j))_{ij}\in F^{n\times n}$ having degree sequence $\mathbf{D}$, i.e. finding \emph{graphical} $\mathbf{D}$. We discuss three classes of graphical degree sequences: those based on the configuration model (CM) and those implied by the Gale-Ryser (GR) and Erd\"{o}s-Gallai (EG) theorems. Given a GR or ER graphical degree sequence $\mathbf{D}$, we can attain $\mathbf{A}$ as the solution to a corresponding maximum-flow problem using for instance the polynomial time Ford-Fulkerson algorithm. 

\subsubsection{CM-graphical sequences}
For $F=\N_{\ge0}$, then one such $\phi$ is given by the \emph{configuration model} (CM) class of random graphs \citep{networks}, which may contain multi-edges or self-edges. Note that the configuration model requires the total degree $D_1+\dotsb+D_n=m$ to be even. The degree sequence $\mathbf{D}=\{D_1,\dotsb,D_n\}$ is said to be \emph{CM-graphical} iff the total degree $D_1+\dotsb+D_n$ is even. Note that CM graphs can be extended to directed graphs.


\subsubsection{GR-graphical sequences}
For $F=\{0,1\}$, the criteria for the existence of $\phi$ (or equivalently, adjacency matrix $\mathbf{A}\in\{0,1\}^{n\times n}$) is given by the Gale-Ryser (GR) theorem in terms of integral vectors, their conjugates, and majorization \citep{gale}. Here, (integral) vector $\mathbf{D}$ is taken to be sorted (descending), and the \emph{conjugate} relation \begin{equation}\label{eq:conj}D_k^* = |\{i: D_i\ge k\}|\for k\in\{1,\dotsb,n\}\end{equation} defines the conjugate (integral) vector $\mathbf{D}^* =\{D_i^*: i=1,\dotsb,n\}$ of $\mathbf{D}$. We say that $\mathbf{D}$ is \emph{majorized} by $\mathbf{D}^*$ if (i) $\sum_{i}^kD_i\le\sum_i^kD_i^*$ for $k=1,\dotsb,n$ and (ii) $\sum_i^nD_i=\sum_i^nD_i^*$. Because the condition $\sum_i^nD_i=\sum_i^nD_i^*$ is automatic, then the existence of $\phi:E\times E\mapsto\{0,1\}$ (or $\mathbf{A}\in\{0,1\}^{n\times n}$) is determined by the truth of condition (i) following from GR. Putting these together, we have the following result. 

\begin{thm}[GR]\label{thm:exist} Consider a non-increasing degree sequence $\mathbf{D}=\{D_i: i=1,\dotsb,n\}$ and its conjugate $\mathbf{D}^* =\{D_i^*: i=1,\dotsb,n\}$ \eqref{eq:conj}. Then a symmetric adjacency matrix $\mathbf{A}\in\{0,1\}^{n\times n}$ with degree sequence $\mathbf{D}$ exists if and only if $\sum_{i}^kD_i\le\sum_i^kD_i^*$ holds for $k=1,\dotsb,n$.
\end{thm} 

In view of this theorem, we say the (non-increasing) degree sequence $\mathbf{D}=\{D_1,\dotsb,D_n\}$ is \emph{GR-graphical} iff $\sum_{i}^kD_i\le\sum_i^kD_i^*$ holds for $k=1,\dotsb,n$. Note that a GR theorem may be formulated to prove the truth of the existence of the adjacency matrix of a directed graph given in-degree and out-degree vectors.

\subsubsection{EG-graphical sequences} Also for $F=\{0,1\}$, we have the following result due to Erd\"{o}s and Gallai \citep{er2}. It gives necessary and sufficient conditions for the existence of a \emph{simple} undirected graph with given degree sequence.

\begin{thm}[EG]\label{thm:exist2} Consider a non-increasing degree sequence $\mathbf{D}=\{D_i: i=1,\dotsb,n\}$. Then a simple symmetric adjacency matrix $\mathbf{A}\in\{0,1\}^{n\times n}$ with degree sequence $\mathbf{D}$ exists if and only if $D_1+\dotsb+D_n$ is even and $\sum_{i}^kD_i\le k(k-1)+\sum_{i=k+1}^n\min(d_i,k)$ holds for $k=1,\dotsb,n$.
\end{thm} 

Hence we say that the (non-increasing) degree sequence $\mathbf{D}=\{D_1,\dotsb,D_n\}$ is \emph{EG-graphical} iff $D_1+\dots+D_n$ is even and $\sum_{i}^kD_i\le k(k-1)+\sum_{i=k+1}^n\min(d_i,k)$ holds for $k=1,\dotsb,n$.


\subsection{Bernoulli thinning of adjacency arrays}\label{sec:thin} Let $E=D=\N_{\ge1}$ and consider the infinite adjacency matrix $\mathbf{A}=(A_{ij})_{ij}\in\R_{\ge0}^{\infty\times\infty}$. Consider the FAIW graph $G$ formed from $(M,f,I)$ with $f\in(\mathscr{E}\otimes\mathscr{E})_{\ge0}$ defined as $f(i,j)=A_{ij}$ having independent Bernoulli weights $W=\{W_x: x\in D\}$, $W_x\sim\kappa_x=\text{Bernoulli}(p_x)$. Then $Mf=W\mathbf{A}W^\intercal=\norm{\diag(W)\mathbf{A}\diag(W)}_1=\norm{\mathbf{A}_W}_1$, where $\mathbf{A}_W=T_W(\mathbf{A})\equiv \diag(W)\mathbf{A}\diag(W)\equiv\{W_xW_yA_{xy}:(x,y)\in D\times D\}$ is the $W$-transform $T_W$ of $\mathbf{A}$ (discussed in Remark~\ref{re:wtransform}). That is, the Bernoulli-thinned adjacency matrix of $G$ is the $W$-transform of $\mathbf{A}$. Essentially, $G$ makes active (or selects) a (possibly infinite) subset of vertices (rows and columns) of $\mathbf{A}$ according to successes of an infinite sequence of independent Bernoulli random variates $\{W_x\}$ with success probabilities $\{p_x\}$. 

Per the example in Remark~\ref{re:finite}, if $p_x=1/x^s$ for $s>1$, then $\E M$ is finite and so $M$ is finite almost surely; hence this parameterization of Bernoulli thinning is finite almost surely. 

\section{Applications}\label{sec:app} We give some applications of the results. First we discuss graphon identification from random unlabeled adjacency matrices (\ref{sec:graphons}), wherein we use the degree distribution as a pseudo-likelihood function of observed degrees. In the second application we formulate prime graphs as STC random graphs with atomic (integer) label distribution supported on the pairs of distinct primes (\ref{sec:primes}). For the third, we define spin networks based on FAIW random graphs and discuss correspondence between the Laplace transform of the graph weight and the partition function of a Gibbs random field (\ref{sec:spin}). In the fourth we formulate random Bayesian networks based on STC random directed acyclic graphs and describe a Markov Chain Monte Carlo algorithm for their inference (\ref{sec:bayes}). Lastly we use STC random graphs to define the architecture of deep neural networks (\ref{sec:neural}).
 
\subsection{Identifying graphons}\label{sec:graphons} Suppose we have a collection of $n$ random realizations $G_1,\dotsb,G_n$ of a random graph $G(\phi)$ model, where $\phi(x,y)\sim\text{Bernoulli}(f(x,y))$ for graphon $f\in(\mathscr{E}\otimes\mathscr{E})_{[0,1]}$. We observe their adjacency matrices $\mathbf{A}_1,\dotsb,\mathbf{A}_n$ but not any information about the vertex labels, where $\mathbf{A}_i\in\{0,1\}^{K_i\times K_i}$ for the $i$. Hence we assume canonical $\nu=\Leb$ with $E=[0,1]$. The goal is to identify the graphon $f$ from the $n$ adjacency matrices. 

Towards identifying the counting distribution $\kappa$, we estimate $c$ and $\delta^2$ from the $\{K_i:i=1,\dotsb,n\}$. The value of $\delta^2-c$ uniquely identifies a Poisson-type (PT) distribution and hence the choice of $\psi$. Then, in possession of PT pgf $\psi$, we use Theorem~\ref{prop:dis} to attain the degree distribution of $Y=d(X)$ for $X\sim\nu$ with pgf \[\psi_Y(t) = \int_E\nu(\D x)\psi(1-\nu f(x,\cdot)+\nu f(x,\cdot) t)\] By Theorem~\ref{thm:dis}, the probability $\P(Y=k)$ may be attained through the Cauchy formula as  \[\P(Y=k)=\frac{1}{2\pi i}\oint_C\frac{\psi_Y(t)}{t^{k+1}}\D t\for k\in\N_{\ge0}\]  

Let $d_{ij}$ be the degree of vertex $j$ of degree vector $\mathbf{D}_i$. Define a \emph{pseudo-log-likelihood} function $\mathcal{L}$ of $f$ given $(\mathbf{D}_1,\dotsb,\mathbf{D}_n)$ as \[\mathcal{L}(f|\mathbf{D}_1,\dotsb,\mathbf{D}_n)=\sum_{i}^n\sum_{j}^{K_i}\log(\P(Y=d_{ij}))\] The condition $f(x,x)=0$, $x\in[0,1]$, is assumed when no self-edges are observed in view of the adjacency matrices. Similarly, $f(x,y)=f(y,x)$ is assumed when all the observed adjacency matrices are symmetric. Then we seek \[\hat{f} = \argmax_{f\in(\mathscr{E}\otimes\mathscr{E})_{[0,1]}}\mathcal{L}(f|\mathbf{D}_1,\dotsb,\mathbf{D}_n)\] 

This approach is based on using the degree distribution to identify the graphon. Because $X$ and $1-X$ have the same distribution for $X\sim\nu=\Leb$, then the integrals of $f(x,y)$, $f(x,1-y)$, $f(1-x,y)$, $f(1-x,1-y)$ with respect to $\nu\times\nu$ are identical. This is essentially relabeling of the vertices by a measure-preserving transformation. This means that identification schemes based on integration create a symmetry in the space of graphons, where there are four solutions in general and two, $\hat{f}(x,y)$ and $\hat{f}(1-x,1-y)$, for symmetric $f$. Theoretically, because there are multiple possible solutions, equally likely, one could say we have learned nothing at this point about the graphon, despite the efforts. However, with some kind of additional information about the graphon, such as monotonicity, e.g., monotone increasing or decreasing, the true solution can be discriminated. 

 Because $f\in L^2([0,1]\times[0,1])$ (it is bounded on a bounded domain), we can expand $f$ as \[f(x,y) =\sum_i\sum_j\beta_{ij}\phi_i(x)\phi_j(y)\for (x,y)\in[0,1]\times[0,1]\] in terms of  the orthonormal Legendre polynomials $\{\phi_i:i\ge1\}$ with $\phi_1(x)=1$ and coefficients $\{\beta_{ij}: i,j\ge1\}$. Consider $m$ bases, organize the coefficients and bases into vectors $\theta=\{\beta_{ij}: i,j=1,\dotsb,m\}\in\Theta=\R^{m^2}$ and $\Phi=\{\phi_{i}(x)\phi_j(y): i,j=1,\dotsb,m\}$ and note $f=\theta\cdot\Phi$. Put $\Theta'=\{\theta\in\Theta: \theta\cdot\Phi\in[0,1]\}$ as the feasible parameter set. Then we solve \[\hat{\theta} = \argmax_{\theta\in\Theta'}\mathcal{L}(f=\theta\cdot\Phi|\mathbf{D}_1,\dotsb,\mathbf{D}_n)\] Note that the first coefficient, $\beta_{11}$, can be identified from the observed degree distributions, where $\beta_{11}=(\nu\times\nu)f=\E Y/c$. Consider $f$ separable (so symmetric): $f(x,y)=g(x)g(y)$. Then we expand $g$ as $g(x)=\sum_i^m\beta_i\phi_i(x)$, so that $f=\theta\cdot\Phi$ where $\theta=\{\beta_{ij}=\beta_i\beta_j:i,j=1,\dotsb,m\}$ and $\Phi=\{\phi_i(x)\phi_j(y): i,j=1,\dotsb,m\}$. Then $\beta_{11}=\beta_{1}^2=\E Y/c$. In estimating $\theta$ through an iterative algorithm, a universal initial condition is $\theta_0 = (\beta_{11},0,\dotsb,0)$.

For example, consider $n=5$ random graphs drawn as $G(\phi)$ where $\phi(x,y)\sim\text{Bernoulli}(f(x,y))$ with $f(x,y)=\ind{}(x\ne y)(1+x)^{-2}(1+y)^{-2}$ and $\kappa=\text{Poisson}(c=30)$. We estimate $c$ from $\{K_i\}=\{31,28,33,25,35\}$ as $c=152/5$ and assume Poisson. We set $\beta_{11} =\E Y/c \simeq 0.513183$ from the mean degrees $\E Y$ and the mean number of vertices across the realizations $c$. We assume that $\kappa$ is Poisson and that $f$ is separable (so symmetric) and monotone decreasing. We use an approximating product (symmetric) quadratic Legendre polynomial system ($m=3$) to approximate $f$; thus $\beta_{11}=\beta_1^2$. We run Metropolis-Hastings on $\theta=(\beta_2,\beta_3)$ for $l=50$ iterations using initial condition $\theta_0=(0,0)$ and proposal $P((x,y),\cdot)=\text{Gaussian}(x,0.01^2)\times\text{Gaussian}(y,0.01^2)$. The trajectory of $\theta$ is $(\theta_i: 0\le i\le l)$. The parameter of the maximum of the maximum likelihoods across the iterates is $\hat{\theta}=(\hat{\beta}_2,\hat{\beta}_3)\simeq(0.2102, 0.0162)$. This solution is monotone increasing, so the true solution in view of assumed $f$ is $\hat{f}(x,y)=(\hat{\theta}\cdot\Phi)(1-x,1-y)$. The relative errors in $L^1$ and $L^2$ are $(\nu\times\nu)|f-\hat{f}| / (\nu\times\nu)f\simeq 0.1122$ and $(\nu\times\nu)(f-\hat{f})^2 / (\nu\times\nu)f^2\simeq 0.0142$, which are reasonably small. Below in Figures~\ref{fig:graphon1} and \ref{fig:graphon2} we show the estimated and true graphon. They are quite similar.

%

\begin{figure}[h!]
\centering
\begingroup
\captionsetup[subfigure]{width=3in,font=normalsize}
\subfloat[Estimated graphon $\hat{f}$\label{fig:graphon1}]{\includegraphics[width=3in]{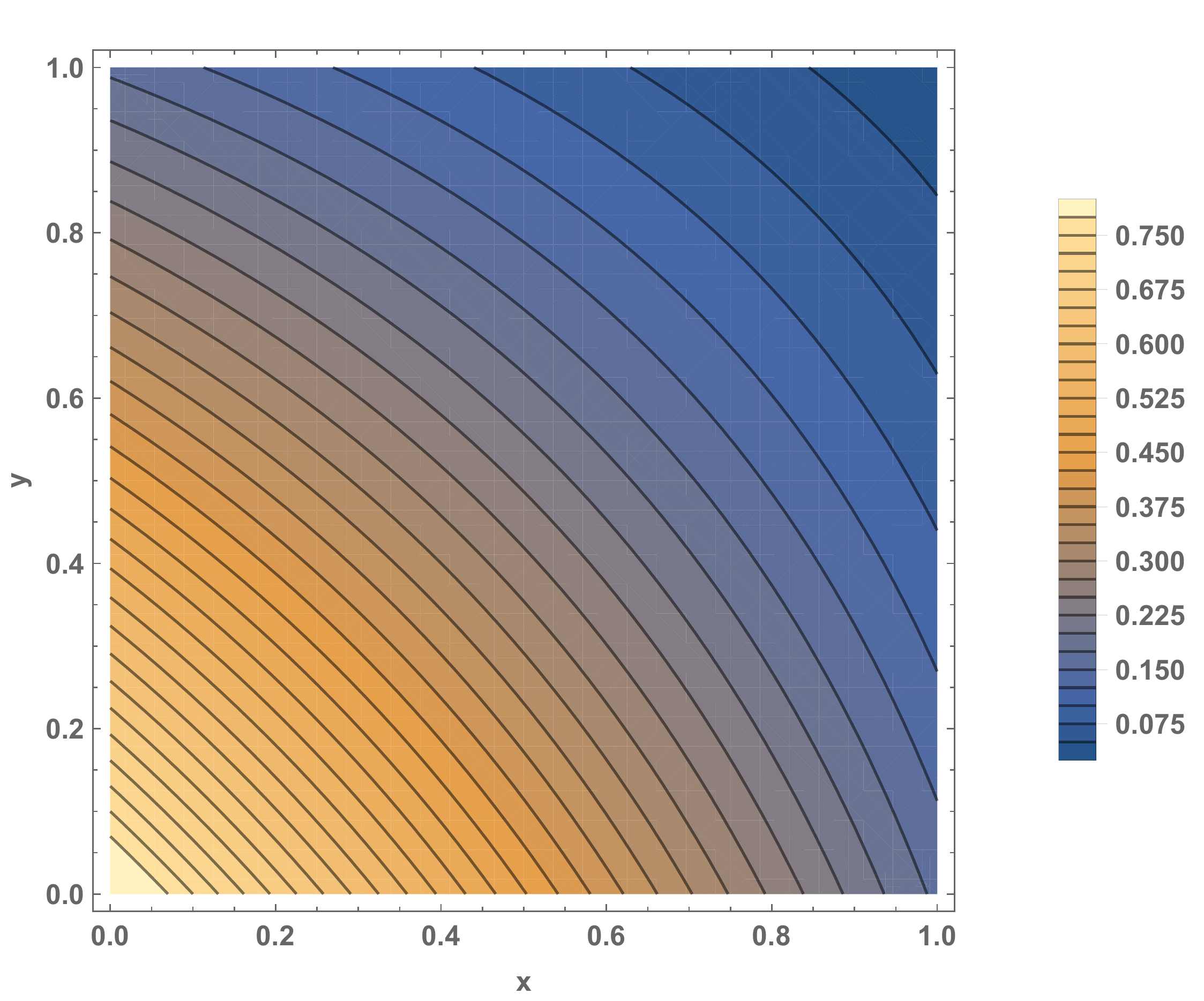}}
\subfloat[True graphon $f$\label{fig:graphon2}]{\includegraphics[width=3in]{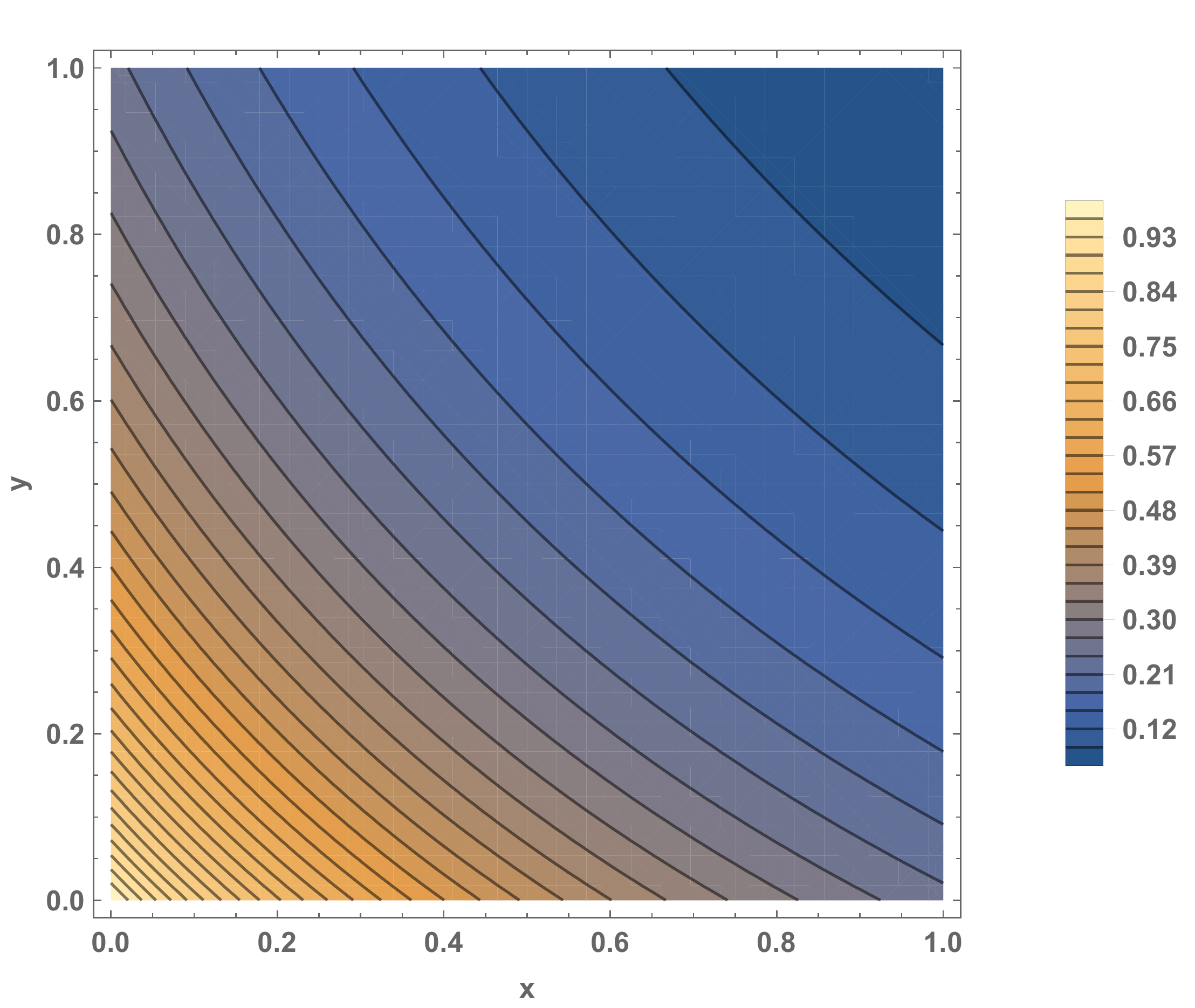}}\\
\endgroup
\caption{Graphon estimation}\label{fig:graphon0}
\end{figure}

\FloatBarrier

Other approaches to graphon estimation including non-parametric estimation of $f$ up to an equivalence class \citep{graphon1,graphon2}, somewhat similar in spirit to the described method. Another approach is to estimate the underlying probability matrix induced by $f$ \citep{graphon3}.

\subsection{Prime graphs}\label{sec:primes} 

Consider $E=\{1,2,\dotsb\}$ and let $P\subset E$ be the collection of the primes. 

Let $\pi$ be the prime counting function, let $\zeta$ be the zeta function $\zeta(s)=\sum_{x\in E}x^{-s}$, $s>1$, and let $\wp$ be the prime zeta function $\wp(s)=\sum_{x\in P}x^{-s}$, $s>1$. 

For zeta distribution $\nu=\text{Zeta}(s)$, $s>1$, we have (zeta) density of the primes $\nu(P)=\wp(s)/\zeta(s)$. $\nu(P)$ is shown below in Figure~\ref{fig:density1}, which has a single maximum around $s\simeq 1.49107$, giving $\nu(P)\simeq0.325236$, and goes to zero as $s\rightarrow\infty$ and as $s\rightarrow1$. For uniform distribution $\nu=\text{Uniform}\{1,\dotsb,n\}$, we have prime density $\nu(P)=\pi(n)/n$, which goes to zero as $n\rightarrow\infty$ as the zero density of the primes and for $n=1$. This can be seen by the prime number theorem (PNT), where $\pi(n)\sim n/\log(n)$.

%

Now consider subset $A \subset E\times E$ as $A=\{(x,y)\in P\times P: x\ne y\}$ so that $\ind{A}(x,y)=\ind{}(x\ne y)\ind{P\times P}(x,y)$ for $(x,y)\in E\times E$. For zeta $\nu$, we have \[(\nu\times\nu)(A) = \frac{\wp^2(s)-\wp(2s)}{\zeta^2(s)}\] and \[\nu \ind{A}(x,\cdot) = \ind{P}(x)\frac{\wp(s)-x^{-s}}{\zeta(s)}\for x\in E\] $(\nu\times\nu)(A)$ is shown below in Figure~\ref{fig:density2} and has a maximum around $s\simeq 1.41152$ of approximately $(\nu\times\nu)(A)\simeq 0.0819344$ and goes to zero as $s\rightarrow\infty$ and as $s\rightarrow1$. For uniform $\nu$ with $n\ge2$, \[(\nu\times\nu)(A) = \frac{\pi(n)(\pi(n)-1)}{n^2}\] and \[\nu \ind{A}(x,\cdot) = \ind{P}(x)(\pi(n)-1)/n\for x\in E\] Note that for uniform $\nu$ both $(\nu\times\nu)(A)\rightarrow0$ and $\nu \ind{A}(x,\cdot)\rightarrow0$ as $n\rightarrow\infty$, again following from PNT.


\begin{figure}[h!]
\centering
\begingroup
\captionsetup[subfigure]{width=3in,font=normalsize}
\subfloat[Prime density $\nu(P)=\wp(s)/\zeta(s)$\label{fig:density1}]{\includegraphics[width=3in]{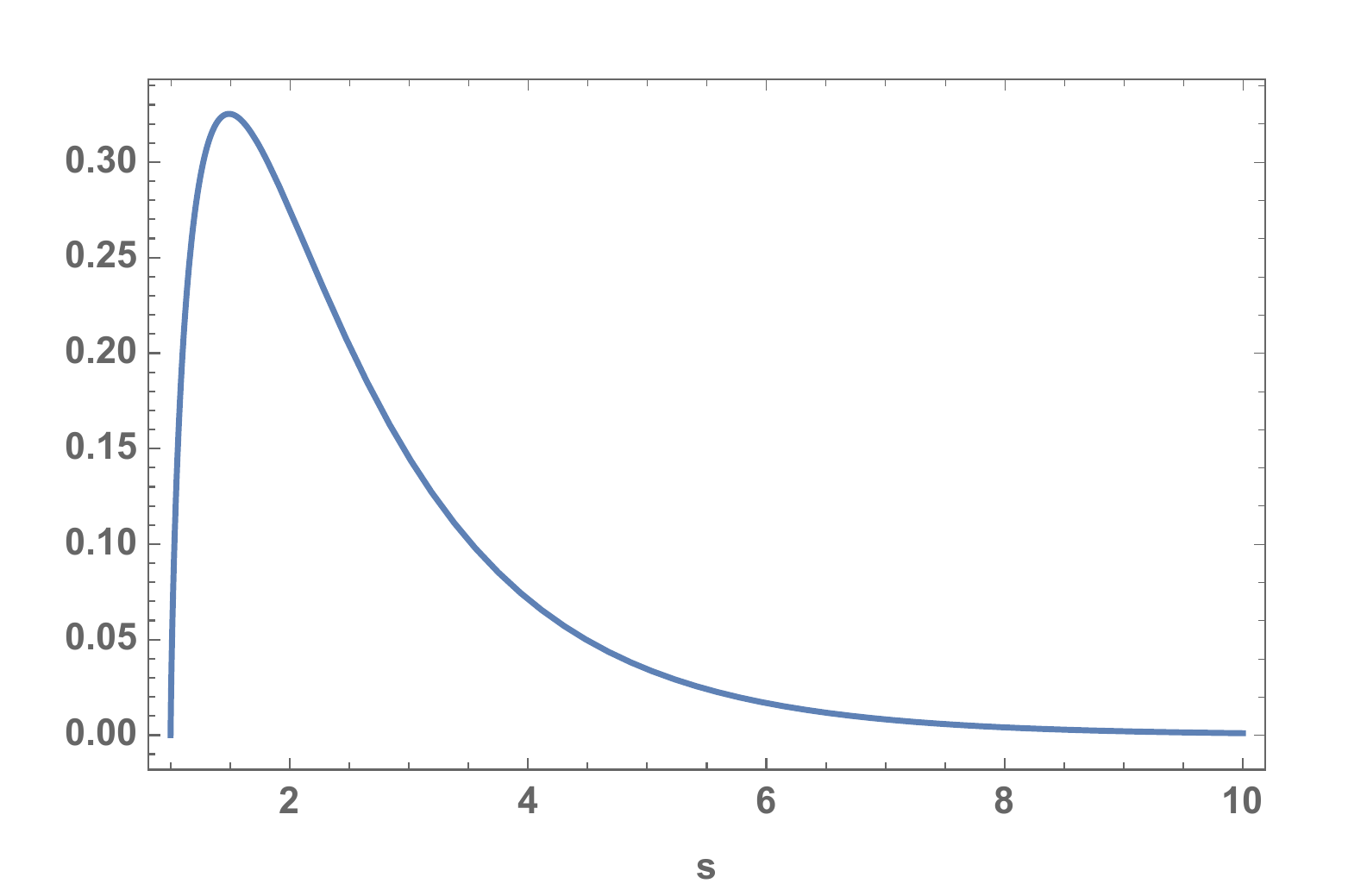}}
\subfloat[Edge density $(\nu\times\nu)(A)=(\wp^2(s)-\wp(2s))/\zeta^2(s)$\label{fig:density2}]{\includegraphics[width=3in]{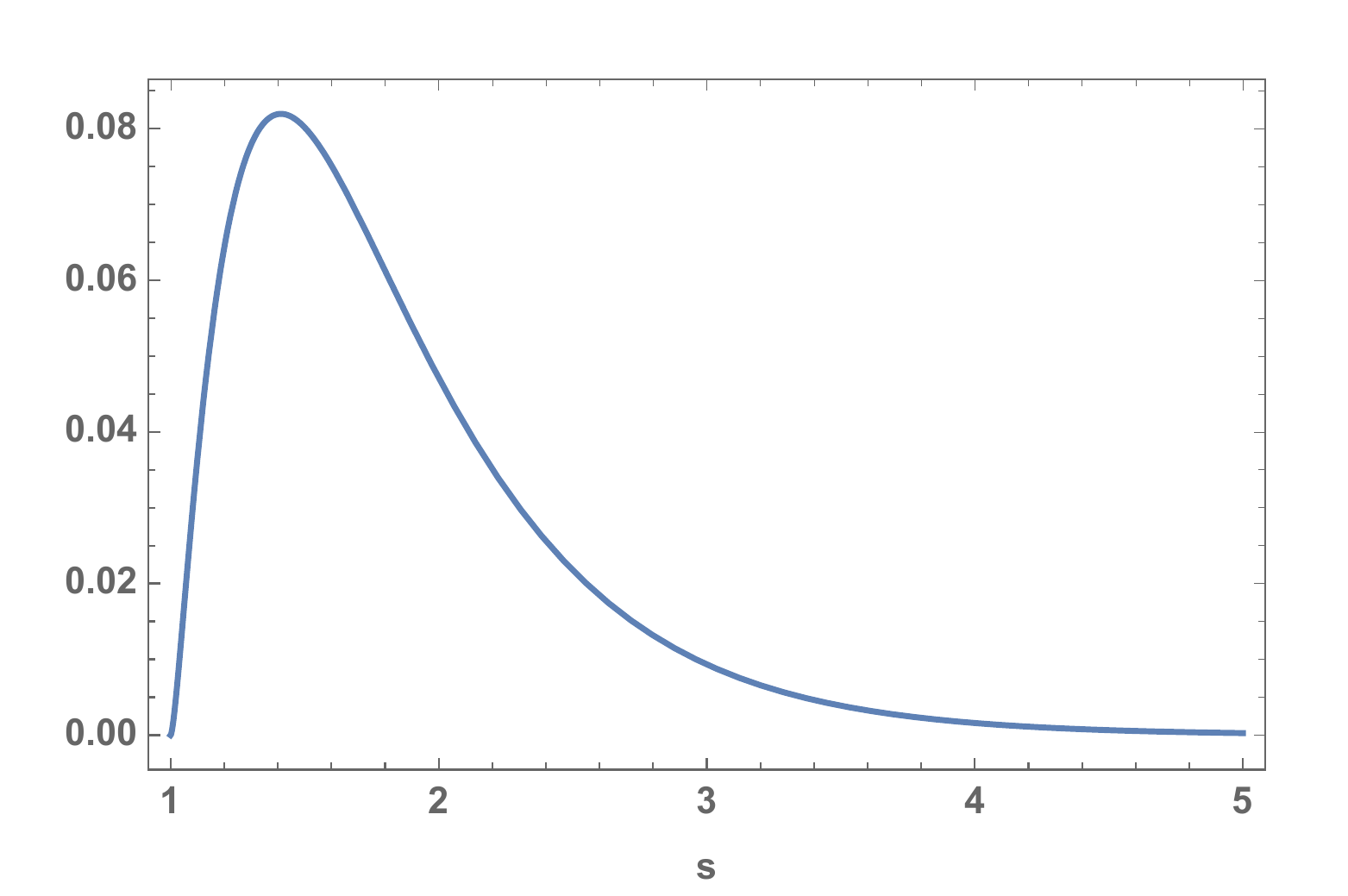}}\\
\endgroup
\caption{Prime and edge densities for zeta $\nu$}\label{fig:densities}
\end{figure}

\FloatBarrier

Now consider the random measure $N=(\kappa,\nu)$ on $(E,\mathscr{E})$ formed by $\mathbf{X}$ and $M=N\times N$ on $(E\times E,\mathscr{E}\otimes\mathscr{E})$ formed by $\mathbf{X}\times\mathbf{X}$. The triple $(M,\ind{A},I)$ forms the random graph $G(\ind{A})=(\mathbf{X},\mathbf{X}\times\mathbf{X},\ind{A})$. The random graph is called \emph{prime} because its support is expressed in terms of the primes. Because the edge space $\mathbf{X}=\{X_i:i=1,\dotsb,K\}$ is a multiset of $K$ iid discrete (labeled) random variables, the edge space is also a multiset where edges can be repeated. The degree function is given by \[d(x) = N\ind{A}(x,\cdot)=K_A^x=\sum_i^K\ind{A}(x,X_i)\for x\in E\] The mean and variance of the degree function are \begin{align*}\E d(x) &= c\nu(\ind{A}(x,\cdot))\\\Var d(x) &= c\nu(\ind{A}(x,\cdot)) + (\delta^2-c)(\nu(\ind{A}(x,\cdot)))^2\end{align*} Put $a(x)=\nu\ind{A}(x,\cdot)$. The pgf of $d(x)$ is \[\psi_A^x(t) = \psi(1-a(x)+a(x)t)\] 

For zeta $\nu$, we have $a(x)=\ind{P}(x)\frac{\wp(s)-x^{-s}}{\zeta(s)}$, so \begin{align*}\E d(x) &= c\ind{P}(x)\frac{\wp(s)-x^{-s}}{\zeta(s)}\\\Var d(x) &= c\ind{P}(x)\frac{\wp(s)-x^{-s}}{\zeta(s)}+(\delta^2-c)\ind{P}(x)\left(\frac{\wp(s)-x^{-s}}{\zeta(s)}\right)^2\end{align*} The number of active vertices $v(G)$ of $G=G(\ind{A})$ for Poisson $\kappa$ and zeta $\nu$ has mean \begin{align*}\E v(G) &= c\nu(1-e^{-ca})\\&=c(1-\sum_{x\in E}\kappa\{x\}e^{-ca(x)})\\&=c\left(1-\sum_{x\in E}\frac{x^{-s}}{\zeta(s)}\exp_-c\ind{P}(x)(\wp(s)-x^{-s})/\zeta(s)\right)\\&=c\left(1-(1-\wp(s)/\zeta(s)+\sum_{x\in P}\frac{x^{-s}}{\zeta(s)}\exp_-c(\wp(s)-x^{-s})/\zeta(s))\right)\\&=c\left(\wp(s)/\zeta(s)-\sum_{x\in P}\frac{x^{-s}}{\zeta(s)}\exp_-c(\wp(s)-x^{-s})/\zeta(s))\right)\\&\sim c\nu(P)\quad\text{as}\quad c\rightarrow\infty\end{align*} and hence is proportional for large graphs to the zeta density of the primes, $\nu(P)$. The same argument can be applied for uniform $\nu$ on $\{1,\dotsb,n\}\subset E$, which gives $\E v(G)\sim c\nu(P)=c\pi(n)/n$ as $n\rightarrow\infty$. 

\begin{figure}[h!]
\centering
\begingroup
\captionsetup[subfigure]{width=3in,font=normalsize}
\subfloat[Zeta\label{fig:gc1}]{\includegraphics[width=3in]{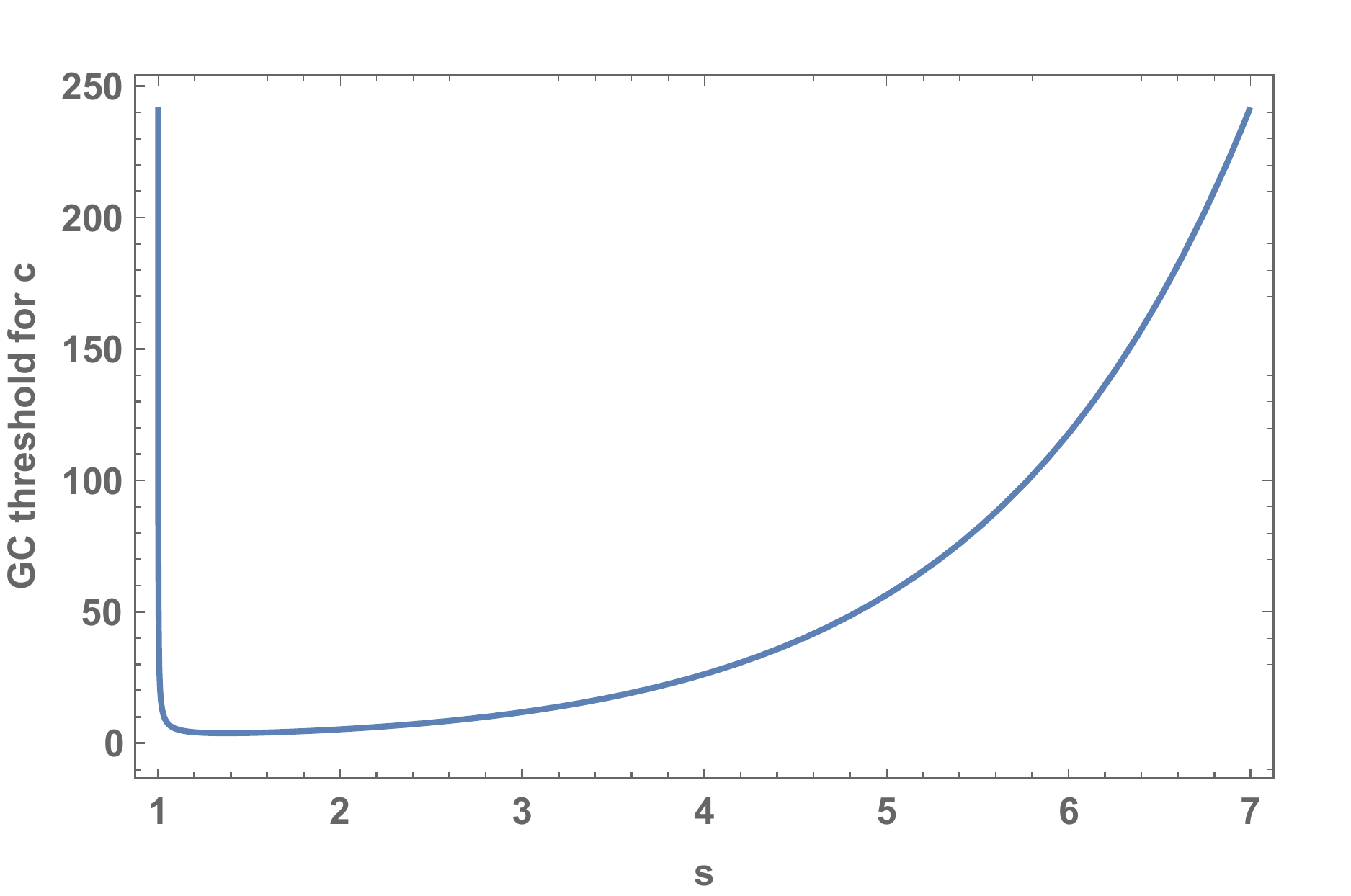}}
\subfloat[Uniform\label{fig:gc2}]{\includegraphics[width=3in]{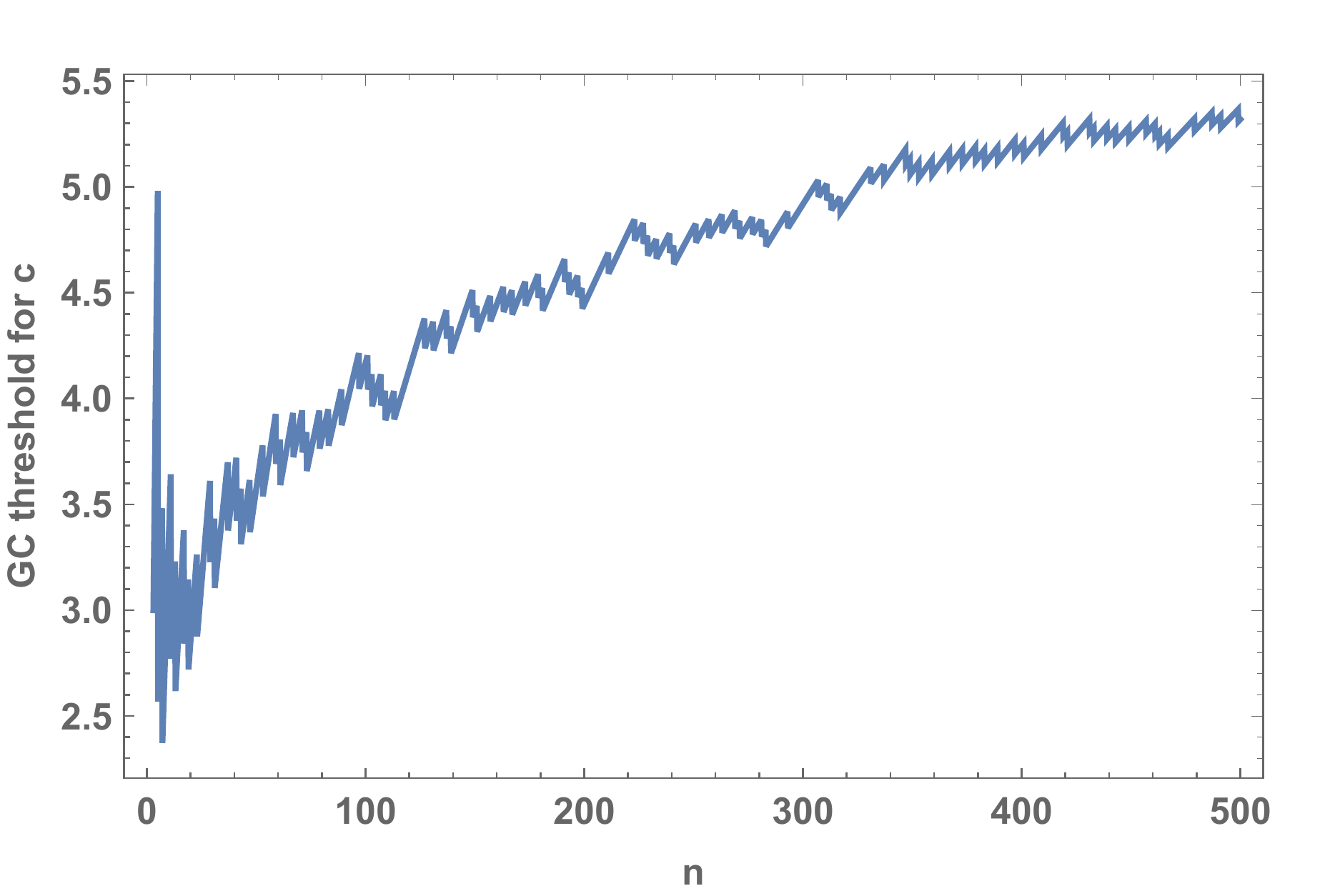}}\\
\endgroup
\caption{Giant component thresholds for zeta and uniform $\nu$ and Poisson $\kappa$}\label{fig:gcs}
\end{figure}

Consider the degree of a random vertex $Y=d(X)$, $X\sim\nu$. For Poisson $\kappa$ and  zeta $\nu$, the mean and variance are \begin{align*} \E Y &= c \frac{\wp^2(s)-\wp(2s)}{\zeta^2(s)}\\\Var Y &= \frac{c}{\zeta^4(s)}\left(c(\wp^3(s)-2\wp(s)\wp(2s)+\wp(3s))\zeta(s) - c(\wp^2(s)-\wp(2s))^2+(\wp^2(s)-\wp(2s))\zeta^2(s)\right)\end{align*} A GC exists almost surely iff \[c > \frac{(\wp^2(s)-\wp(2s))\zeta(s)}{\wp^3(s)-2\wp(s)\wp(2s)+\wp(3s)}\] For Poisson $\kappa$ and uniform $\nu$, the GC condition is \[n\ge 3,\quad c > n/(\pi(n)-1)\sim\log(n)\text{ as }n\rightarrow\infty\] In Figures~\ref{fig:gc1} and \ref{fig:gc2} we show the giant component thresholds for $c$ for zeta and uniform $\nu$. The zeta result in Figure~\ref{fig:gc1} shows that for both $s\rightarrow1$ and $s\rightarrow\infty$, the mean number of vertices must increase $c\rightarrow\infty$ in order to guarantee the existence of a GC. In particular the value of $c$ rapidly (exponentially) grows in $s$. For uniform $\nu$, Figure~\ref{fig:gc2} shows that GC threshold drops precipitously and then slowly (logarithmically) increases. 

%

\FloatBarrier

\subsection{Spin networks}\label{sec:spin} We illustrate that a FAIW random graph forms a Gibbs random field of non-negative integer-valued spins on a lattice, a spin network.

Consider a FAIW product random measure $M$ on $(E\times E,\mathscr{E}\otimes\mathscr{E})$ where $E=\N_{\ge1}^d$ is a $d$-dimensional \emph{lattice}. Without loss of generality, we take $d=1$. Let $D=E$ be its atoms (the lattice) and let $W=\{W_x: x\in D\}$ be independent non-negative integer weights (\emph{spins}) in $H=\N_{\ge0}^{|D|}$ distributed $W_x\sim\kappa_x$ with mean $c_x$ and variance $\delta_x^2$. We refer to $W$ as a \emph{spin configuration}, $H$ the \emph{spin configuration space}, and the zero weight as the \emph{ground state} of a spin. Take subset $A\subseteq E\times E$, such as through \emph{$n$-nearest-neighbors} as $A=\{(x,y)\in E\times E: 0<\norm{x-y}\le n\}$. Let deterministic $f:E\times E\mapsto\R_{\ge0}$ represent an \emph{interaction} function, where $f(x,y)$ is the interaction of locations $x$ and $y$. Define \emph{local interaction} function $g = f\ind{A}$ and weight function $h_W$ as $h_W(x,y)=W_{x}W_{y}g(x,y)$. The random graph $G_W=G(h_W) = (D,D\times D, h_W)$ is called a ($A$-local) \emph{spin network} with spin configuration $W$. The random variable $E_W=Mg$ formed as \begin{align*}E_W = Mg &= \sum_{(x,y)\in D\times D}W_xW_yf(x,y)\ind{A}(x,y)\\&=\sum_{(x,y)\in (D\times D)\cap A}W_xW_yf(x,y)\end{align*} is the \emph{energy} (total interaction) of the ($A$-local) spin network $G_W$, having mean \[\E E_W = \E Mg = \sum_{(x,y)\in (D\times D)\cap A}Z_{xy}f(x,y)\] where \[Z_{xy} = \begin{cases}c^2_x+\delta_x^2&\text{if }x=y\\c_xc_y&\text{otherwise}\end{cases}\] 

Let $B=g(D\times D)$ be the local potential double array, where $Mf = WBW^{\intercal}$. The $W$-transform $T_W$ of $B$ is the spin-weighted local potential double array $B_W=T_W(B) = \diag(W)B\diag(W)=\{W_xW_yB_{xy}: (x,y)\in D\times D\}$. This is the adjacency array of the ($A$-local) spin network $G_W$. 

Following from Theorem~\ref{prop:lf}, for general $h\in(\mathscr{E}\otimes\mathscr{E})_{\ge0}$, the Laplace transform of $Mh$ is given as \[\varphi(\alpha) = \E e^{-\alpha Mh} = \sum_{w\in H}\P(W=w)e^{-\alpha E_w}=\sum_{w\in H}\left(\prod_{x\in D}\kappa_x\{w_x\}\right)e^{-\alpha\sum_{(x,y)\in D\times D}w_xw_yh(x,y)}\for \alpha\in\R_{\ge0}\] Per Remark~\ref{re:part}, $\varphi$ is the partition function $\mathcal{Z}$ of a Gibbs random field on a lattice. Therefore, the FAIW random graph construction of spin networks is identical to those generated by the Gibbs random field on a lattice containing non-negative integer-valued spins and no external interactions. 

Now we consider a slightly more general spin network: one containing external interactions. We use the generalized product random measure  of Proposition~\ref{prop:gen} and the Laplace functional. Consider interaction function $h\in(\mathscr{E}\otimes\mathscr{E})_{\ge0}$ and energy $E_W=Mh$. Let $k\in\mathscr{E}_{\ge0}$ be an \emph{external interaction} function with energy $E_W^*=Nk$ and define the overall energy $E^+_W = E_W + E_W^*$. The overall energy has mean $\E E_W^+ = \E E_W + \E E_W^*$ and Laplace transform \[\varphi^+(\alpha) = \E e^{-\alpha E_W^+}=\sum_{w\in H}\left(\prod_{x\in D}\kappa_x\{w_x\}\right)e^{-\alpha\left(\sum_{(x,y)\in D\times D}w_xw_yh(x,y) + \sum_{x\in D}w_xk(x)\right)}\for \alpha\in\R_{\ge0}\] As before, $\varphi^+$ is the partition function $\mathcal{Z}^+$ of a Gibbs measure with energy function $E_W^+$.

\subsection{Bayesian networks}\label{sec:bayes} A Bayesian network (BN) is probabilistic graphical model expressed in terms of a directed acyclic graph (DAG) whose vertices index random variables and whose edges index conditional dependencies. We assume the DAG is finite with $n$ vertices. Let $X=(X_1,\dotsb,X_n)$ be a vector of $n$ random variables taking values in set $F$ of measurable space $(F,\mathscr{F})$. Let $\text{pa}(v)$ be the set of parent vertices of vertex $v$, and let $X_{\text{pa}(v)}=(X_i: i\in\text{pa}(v))$ be the vector of parent random variables corresponding to vertex $v$. The law of $X$ is specified by the BN of DAG $G$ as \[\P_G(X=x)=\prod_{i\in V}\P(X_i=x_i|X_{\text{pa}(i)}=x_{\text{pa}(i)})\for x\in F\] 


We construct random BNs using STC random DAGs. Consider STC random measure $N=(\kappa,\nu)$ on $(E,\mathscr{E})=([0,1],\mathscr{B}_{[0,1]})$ formed by $\mathbf{Z}=\{Z_i:i=1,\dotsb,K\}$, where $Z_i\sim\nu=\Leb$ and $K\sim\kappa$. Consider product STC random measure $M=N\times N$ formed by $\mathbf{Z}\times\mathbf{Z}$. Take subset $A=\{(x,y)\in E\times E: x<y\}$.  Let $G=G(\phi)=(\mathbf{Z},\mathbf{Z}\times\mathbf{Z},\phi)$ be a STC random DAG, where $\phi$ is the Bernoulli transformation with local graphon $f_A=f\ind{A}$ of graphon $f$. Let $(F,\mathscr{F})$ be a measurable space and consider the BN of the random vector $X=(X_z: z\in\mathbf{Z})$ supported on $F$ with law relative to the DAG $G$ \begin{equation}\label{eq:pg}\P_G(X=x)=\prod_{z\in\mathbf{Z}}\P(X_z=x_z|X_{\text{pa}(z)}=x_{\text{pa}(z)})\for x\in F\end{equation} A common setting is to take $\kappa=\text{Dirac}(n)$ and fixed ($n$-dimensional) space $(F,\mathscr{F})$ so that the resulting BNs are relative to common $n$-dimensional $X$ supported on $F$. Note that \[\P(X_z=x_z|X_{\text{pa}(z)}=x_{\text{pa}(z)}) = Q_{\text{pa}(z)}(x_{\text{pa}(z)},x_z)\for x_{\text{pa}(z)}\in F_{\text{pa}(z)},\quad x_z\in F_z \] forms the probability transition kernel $Q_{\text{pa}(z)}:F_{\text{pa}(z)}\mapsto F_z$, so the law $\P_G$ is encoded by the collection of transition kernels $\{Q_{\text{pa}(z)}: z\in\mathbf{Z}\}$. 

Now we discuss an inference scheme of $G$ given observed data based on Metropolis-Hastings (MH) algorithm. Suppose we have an independency $\mathbf{X}$ of $n$-dimensional random vectors in $F$. Define the \emph{likelihood function} $\mathcal{L}$ of $G$ given $\mathbf{X}$ as \[\mathcal{L}(G|\mathbf{X})=\prod_{x\in\mathbf{X}}\P_G(X=x)\] Next define the \emph{proposal kernel} $\pi$ for STC random graphs through the \emph{rewiring transform} $S_n$ \ref{eq:sn} as \[S_n(G)\sim\pi(G,\cdot)\] where $n\in\{1,\dotsb,K\}$ is the number of rewired vertices. Recall that $S_n$ preserves the mean adjacency matrix of $G$. Equipped with an ability to compute (or estimate) the conditional probabilities composing the law $\P_G$ and the likelihood $\mathcal{L}$ and the proposal $\pi$, the MH algorithm generates a sequence of rewired graphs that converges to the posterior distribution of $G$, which sequence preserves the mean adjacency matrix of $G$ for fixed graphon $f$. Variations on the rewiring transform, such as rewiring individual (random) edges, offer flexibility in optimizing the rate of convergence of the MH sampler. The graphon controls the properties of the random DAGs and thus behaves as a functional `prior' on the space of DAGs.

A simple counting representation of the transition kernel $Q_{\text{pa}(z)}$ for $z\in\mathbf{Z}$ is given as follows. Let $D_z$ be the empirical distribution function of $(\mathbf{X}_{\text{pa}(z)},\mathbf{X}_z)$ on $(F_{\text{pa}(z)}\times F_z,\mathscr{F}_{\text{pa}(z)}\otimes\mathscr{F}_z)$. Now let $A_1,\dotsb,A_q$ be a partition of $F_{\text{pa}(z)}$ with $D_{z}(A_i\times F_z)>0$, $i=1,\dotsb,q$, and let $B_1,\dotsb,B_r$ be a partition of $F_z$. Then the transition kernel is formed as \[Q_{\text{pa}(z)}(x_{\text{pa}(z)},x_z)=\sum_{i}^q\sum_{j}^r\frac{D_z(A_i\times B_j)}{D_{z}(A_i\times F_z)}\ind{A_i\times B_j}(x_{\text{pa}(z)},x_z)\] and takes $qr$ possible values. The partition $A_1,\dotsb,A_q$ may be efficiently attained for example using K-means clustering ($K=q$) of $\mathbf{X}_{\text{pa}(z)}$. 

\subsection{Deep neural networks}\label{sec:neural} A key issue in using deep (multi-layer) neural networks (DNNs) is defining their architectures. Classical feed-forward neural networks fully connect neurons between successive layers, where the number of connections from layer $i$ consisting of $n_i$ neurons to any (and every) neuron of layer $i+1$ is $2^{n_i}$. Instead of fully connecting neurons, a strategy is to sparsely (and randomly) connect them, which has been shown to match performance of manually-optimized architectures and to improve computational efficiency \citep{nn}. We sketch a simple stochastic block wiring algorithm for multi-layer feed-forward architectures using STC Bernoulli graphs. We compute the number of total connections, the neuron degrees, and number of parameters. 

Consider an architecture having $n+2$ layers, with one input layer, $n$ hidden layers, and one output layer. We assume the nodes of the input layer are fully connected to each neuron of the first hidden layer, and the neurons of the final hidden layer are fully connected to each of the nodes of the output layer. Hence we focus on random connections among the neurons of the hidden layers. 

Recall the set-up of the STC Bernoulli random graph $G=(\mathbf{X},\mathbf{X}\times\mathbf{X},\phi)$ identified to $(M,\phi,I)$, where $(E,\mathscr{E})$ is a measurable space, $K\sim\kappa$ is a $\N_{\ge0}$-valued random variable, $\mathbf{X}=\{X_i:i=1,\dotsb,K\}$ is an independency $E$-valued random variables with law $\nu$ that forms the random measure $N=(\kappa,\nu)$ on $(E,\mathscr{E})$, $M=N\times N$ is the product random measure on $(E\times E,\mathscr{E}\otimes\mathscr{E})$, and $\phi$ is a Bernoulli transformation with mean $f\in(\mathscr{E}\otimes\mathscr{E})_{[0,1]}$.

Let $(E,\mathscr{E})=(\{1,\dotsb,n\},2^{\{1,\dotsb,n\}})$ be the layer space of neurons with distribution $\nu$. Define Bernoulli transform mean $f$ as $f(x,y)=\ind{}(y=x+1)p(x,y)$ where $p\in(\mathscr{E}\otimes\mathscr{E})_{[0,1]}$. This is a stochastic block model of feed-forward neuron connectivity, where each layer connects to the next. The total number of connections/edges is $e(G) = (M\circ\phi^{-1})I$ with mean \[\E e(G) = (c^2+\delta^2-c)(\nu\times\nu)f = (c^2+\delta^2-c)\sum_{x}p(x,x+1)\nu\{x\}\nu\{x+1\}\] The out- and in-degree functions of the neurons by layer are $d_o(x)=Nf(x,\cdot)$ and $d_i(x)=Nf(\cdot,x)$ with means \begin{align*}\E d_o(x) &= c\nu f(x,\cdot) = c\sum_yf(x,y)\nu\{y\}=cp(x,x+1)\nu\{x+1\}\ind{}(x<n)\for x\in E\\\E d_i(x) &= c\nu f(\cdot,x) = c\sum_yf(y,x)\nu\{y\}=cp(x-1,x)\nu\{x-1\}\ind{}(x>1)\for x\in E\end{align*} For Poisson $\kappa$, the degrees are Poisson random variables. 

Let $F_x$ be a neuron indexed by $x\in\mathbf{X}$. It has input vector $z\in\R^{n_x}$ and is represented as $F_x(z)=\sigma(a_x\cdot z + b_x)$ where $\sigma$ is a sigmoid function, $a_x\in\R^{n_x}$ is a weight vector, and $b_x\in\R$ is a bias term. $F_x$ has $n_x+1$ parameters. Thus the number of parameters corresponding to connections \emph{among} the hidden layers for a feed-forward model is given by $e(G)+K_{>1}$, where $e(G)$ is the number of weight parameters, each edge contributing one weight, and $K_{>1}=N\ind{>1}$ is the number of bias parameters corresponding to the neurons in the second and following layers.



\section{Discussions and conclusions}\label{sec:discuss}



We use the product random counting measures and a (possibly random) transformation to encode a general model for random graphs. This model constructs, through random transformations, directed and undirected random graphs, with or without loops and/or multiedges. Undirected are retrieved using for instance graphon-based Bernoulli transformations. Directed are retrieved using digraphon-based or asymmetric transformations. Multigraphs are generated using the Poisson transformation. A number results are developed, including mean number of edges and edge weight, degree functions, degree distributions, number of active vertices, and so on. 

We apply the results to graphon identification of observed unlabeled graphs, retrieving an estimation scheme under assumptions. Graphon identification is a subtle problem and there are numerous approaches. The method for graphon identification we describe is reasonably general and scalable. 

We study prime graphs as those random graphs with supports expressed in term of the primes $P$. This graph is formed by an indicator function on a discrete space $E\times E$ where $E=\{1,2,\dotsb\}$.  The mean number of active vertices scales with the density of the primes $\nu(P)$ and the criterion for a giant components scales based on the parameter of the underlying discrete distribution $\nu$, either $s>1$ for zeta or $n\ge3$ for uniform.

We construct spin networks from product FAIW random measures. We show the Laplace transform of FAIW random graph network energy is the partition function of a Gibbs random field,  conveying stochastic equivalence of the models. The FAIW mean measure is used to compute the mean energy. 

We formulate random Bayesian networks (BNs) through random directed acyclic graphs (DAGs) generated by product STC random measures. We sketch a Markov Chain Monte Carlo inference algorithm to sample the posterior distribution of the DAG given observed data. 

In the last application we define neural network architectures through product STC random measures using a stochastic block or community model of neuron connectivity. The number of edges/connections controls the total number of network parameters, and the degrees reveal neuron connectivity by layer. 

The class of random 2-graphs formed by product random measures $M=N\times N$ can represent arbitrary 2-graphs, where the mean product measure $\E M$, possibly infinite, encodes the mean edge count and weight, and where the mean measure $\E N$ encodes the mean in- and out-degree functions. Moreover these encodings use only the marginal transition kernel of the transformation, not its law. 

A limitation of this formulation is that the mean measures of higher-order product random measures become unwieldy. So too their Laplace functionals. If instead we take factorial powers of random measures, removing contributions of repeats, then the ensuing mean measures are factorial moment measures, which are easier to work with, albeit with loss of generality. 

We have explored some of the theory of random graphs as formed by product random measures. Given the generic nature of the representation, there are other transforms to explore and extensions to consider, e.g., higher-order (tensor) products of random measures and their superpositions conveying graphs having edges of arbitrary sizes. Another area of future work is developing a limit theory for the STC random graph models of this article.  

\section*{Funding Sources} H.R. acknowledges support from Army Research Office Grant W911NF-19-1-0382. 

 \bibliographystyle{apalike}
\bibliography{../stones3}

\end{document}